\documentclass{article}

\usepackage{graphics,latexsym}
\usepackage{amsmath, amssymb,amsthm}
\usepackage{graphicx}
\usepackage{amsmath,amssymb}
\usepackage{graphics}
\usepackage{graphicx}
\usepackage[latin1]{inputenc}
\usepackage{epsfig}
\usepackage{comment}

\newcommand{\R}{\mathbb{R}}

\newcommand{\Xnewnm}{\tilde x_{nm} }
\newcommand{\Xoldnm}{x_{n}}
\newcommand{\Vnewnm}{\tilde v_{nm} }
\newcommand{\Voldnm}{v_{m}}

\newcommand{\incrXnm}{ v_{m} \, \Delta t}

\newcommand{\incrVnmb}{E_N(x_{n}) \, \Delta t}

\newcommand{\EOD}{\end{document}}

\newcommand{\dv}{dv}

\newcommand{\Xv}{\mathbf{X}}
\newcommand{\Yv}{\mathbf{Y}}

\newcommand{\as}{a}
\newcommand{\bs}{b}
\newcommand{\cs}{c}
\newcommand{\ds}{d}

\newcommand{\fs}{f}
\newcommand{\gs}{g}

\newcommand{\vs}{v}

\newcommand{\xs}{x}

\newcommand{\Bs}{B}

\newcommand{\Es}{E}

\newcommand{\Ps}{P}
\newcommand{\Qs}{Q}

\newcommand{\Vs}{V}

\newcommand{\Xs}{X}

\newcommand{\abs}[1]{\vert#1\vert}

\usepackage{xcolor}

\begin{document}

\begin{Large}
  \begin{center}
    A Semi-Lagrangian Spectral Method\\ \smallskip
    for the Vlasov-Poisson System\\ \smallskip
    based on Fourier, Legendre and Hermite Polynomials
  \end{center}
\end{Large}

\begin{center}
  Lorella Fatone\\
  {\small Dipartimento di Matematica, Universit\`a degli Studi di
    Camerino, Italy,}\\
  {\small lorella.fatone@unicam.it}
\end{center}

\begin{center}
  Daniele Funaro\\
  {\small Dipartimento di Scienze Fisiche, Informatiche e Matematiche,
    Universit\`a degli Studi di Modena e Reggio Emilia, Italy,}
  {\small daniele.funaro@unimore.it}
\end{center}

\begin{center}	
  Gianmarco Manzini\\
  {\small Group T-5, Applied Mathematics and Plasma Physics,
    Theoretical Division, Los Alamos National Laboratory, Los Alamos,
    NM, USA,}
  {\small gmanzini@lanl.gov}
\end{center}

\bigskip

\begin{abstract}
In this work, we apply a semi-Lagrangian spectral method for the Vlasov-Poisson system, previously 
designed for periodic Fourier discretizations, by implementing 
Legendre polynomials and Hermite functions in the approximation of the distribution function with 
respect to the velocity variable. We discuss second-order accurate-in-time schemes, 
obtained by coupling spectral techniques in the space-velocity domain with a BDF 
time-stepping scheme. The resulting method possesses good conservation properties, which 
have been assessed by a series of numerical tests conducted on the standard two-stream 
instability benchmark problem. In the Hermite case, we also investigate the numerical 
behavior in dependence of a scaling parameter in the Gaussian weight. Confirming previous results from 
the literature, our experiments for different representative values of this parameter, indicate 
that a proper choice may significantly impact on accuracy, thus suggesting that suitable strategies
should be developed to automatically update the parameter during the time-advancing procedure.
\end{abstract}

\bigskip

\section{Introduction}

%% introduction
A semi-Lagrangian spectral method has been proposed
in~\cite{Fatone-Funaro-Manzini:2018} for the numerical approximation
of the nonrelativistic Vlasov-Poisson equations.
These equations describe the dynamics of a collisionless plasma of
charged particles, under the effect of a self-consistent electrostatic
field~\cite{Boyd-Sandserson:2003}.
For the exposition's sake, we assume that each plasma species
(electrons and ions) is described by a 1D-1V distribution function,
i.e., a function defined in a phase space that is one-dimensional in
the independent variables $x$ (space) and $v$ (velocity).
The approximation introduced in~\cite{Fatone-Funaro-Manzini:2018} has
been specifically tested on Fourier-Fourier periodic discretizations
for both variables in the phase space.
The main goal of this paper is to extend our approach to
Fourier-Legendre and Fourier-Hermite discretizations, i.e., by
considering the spectral type representation with respect to $v$
provided by Hermite and Legendre polynomials, in conjunction
with a high-order semi-Lagrangian technique in time.

\smallskip

%% 4 families of methods
Semi-Lagrangian methods have been proposed in
different frameworks such as the Finite Volume
method~\cite{Filbet:2001,Banks-Hittinger:2010}, Discontinuous
Galerkin
method~\cite{Ayuso-Carrillo-Shu:2011,Ayuso-Carrillo-Shu:2012,Heath-Gamba-Morrison-Michler:2012},
finite difference methods based on ENO and WENO polynomial
reconstructions~\cite{Christlieb-Guo-Morton-Qiu:2014}, as well as in the  propagation
of solutions along the characteristics in an operator splitting
context~\cite{Arber-Vann:2002,%
  Carrillo-Vecil:2007,%
  Cheng-Knorr:1976,%
  Filbet-Sonnendrucker-Bertrand:2001,%
  Filbet-Sonnendrucker:2003,%
  Sonnendrucker-Roche-Bertrand-Ghizzo:1999,%
  Crouseilles-Respaud-Sonnendrucker:2009}.
% 
% Semi-Lagrangian methods were first developed for meteorological
% applications in the early '90s
% \cite{Bermejo:1991,Bermejo:1995,Staniforth-Cote:1991}.
These methods offer an alternative to Particle-in-Cell (PIC)
methods~\cite{Cottet-Raviart:1984,Wollman:2000,Wollman-Ozizmir:1996,Brackbill:2016,Chen-Chacon:2015,Chen-Chacon-Barnes:2011,Lapenta:2017,Lapenta-2011,Markidis-Lapenta:2011,Taitano-Knoll-Chacon-Chen:2013}
and to the so called Transform methods based on spectral
approximations~\cite{Manzini:2016:JCP:journal,Parker-Dellar:2015,Klimas:1983,Vencels-Delzanno-BoPeng-Laure-Markidis:2015,Vencels:2016:JPCS:journal}.
PIC methods are very popular in the plasma physics community and are
the most widely used methods because of their robustness and relative
simplicity~\cite{Birdsall-Langdon:2005}.
The PIC method has been successfully used to simulate the behavior of
collisionless laboratory and space plasmas and provides excellent
results for the modeling of large scale phenomena in one, two or three
space dimensions~\cite{Birdsall-Langdon:2005}.
%, and implicit and energy
%preserving PIC formulations that are suitable to long time integration
%problems are available.
%%
On the other hand, spectral Transform methods use Hermite basis functions for unbounded
domains, Legendre basis functions for bounded domains, and Fourier
basis functions for periodic domains, and can outperform PIC in
Vlasov-Poisson
benchmarks~\cite{Camporeale-Delzanno-Vergen-Moulton:2015,Camporeale-Delzanno-Lapenta-Daughton:2006}.
Moreover, they can be extended in an almost straightforward way to
multidimensional simulations of more complex models, like
Vlasov-Maxwell~\cite{Delzanno:2015}.
Convergence of various formulations of these methods was shown
in~\cite{Gajewski-Zacharias:1977,Manzini-Funaro-Delzanno:2017:journal}.

\smallskip

The method proposed in~\cite{Fatone-Funaro-Manzini:2018} works as follows.
At each grid point, the updated value is set up to be equal to the 
value obtained by going backward, by a suitably small amount,  
along the local characteristic lines (approximated in a certain way).
The algorithm is the result of a Taylor development of arbitrary accuracy.
This requires the computation of derivatives in the variable $x$ and $v$
of appropriate order. This operation can be carried out at a spectral
convergence rate, therefore the final algorithm may turn out extremely accurate,
depending on how many terms are included in the Taylor development.
%%
%The second order of accuracy in time can be achieved according to the
%Method-of-Lines (MOL) by coupling the spectral phase-space
%discretization with a second-order accurate Backward Differentiation
%Formula (BDF).

%As mentioned above, we extend the Fourier-Fourier Semi-Lagrangian method
%introduced in~\cite{Fatone-Funaro-Manzini:2018} to the case of a
%Fourier discretization in space and Legendre or Hermite discretizations in
%velocity.
%%
%Since a spectral collocation approach for both the space and velocity
%representations, this method provide the spectral accuracy of the
%aforementioned Transform methods.
%%

\smallskip

The literature on spectral methods is rather rich. Theoretical studies
and numerical results concerning applications in the field of partial
differential equations in unbounded domains have been subject of
investigation in the last 30 years.
In this context, functional spaces generated by Laguerre or Hermite
polynomials provide the main framework.
As far as Hermite approximations are concerned, a non exhaustive list
of references is:
\cite{Funaro-Kavian:1990,Guo:1999,Guo-Xu:2000,Tang:1993,Guo-Shen-Xu:2003}.

%%{\textcolor{red}{Qui bisognerebbe aggiungere le references che avevo segnalato. Non le ho trovate nel bib.
%%vedi Blocco 1}}  

\smallskip

The first attempt in using Hermite polynomials to solve the Vlasov
equation dates back to the paper~\cite{Grad:1949}.
In that work, a Hermite basis is used in the velocity variable for
the distribution function of a plasma in a physical state near the
thermodynamic equilibrium, i.e., the Maxwellian (Gaussian) distribution function.
Since Hermite polynomials are orthogonal with respect to the 
exponential weight $\exp(-v^2)$, a close link  to
Maxwellian distribution functions is soon obtained, simplifying
enormously the formulation of the method.
Within this approach, exact conservation laws in the
discrete setting, i.e., discrete invariants in time for number of
particles (mass, charge), momentum and energy, can be constructed from
Hermite expansion's coefficients.
Moreover, as pointed out in~\cite{Vencels-Delzanno-BoPeng-Laure-Markidis:2015,%
Vencels:2016:JPCS:journal,%
Delzanno:2015,Manzini:2016:JCP:journal,%
Manzini-Funaro-Delzanno:2017:journal}, 
expansions in Hermite basis functions
are intrinsically multiscale, providing a natural connection between
low-order moments of the plasma distribution function and typical
fluid moments.
%% so that multiscale effects as in modeling kinetic physics are
%% naturally taken into account.

\smallskip

The weight function of the Hermite basis can be easily generalized by
introducing a parameter $\alpha$ in such a way that it becomes
$\exp(-\alpha^2 v^2)$.
A proper choice of this parameter can significantly improve the
convergence properties of Hermite polynomial
series~\cite{Boyd:1980,Boyd:1984,Tang:1993}.
This fact was also confirmed in earlier works on plasmas physics based
on Hermite spectral methods
(see~\cite{Schumer-Holloway:1998,Holloway:1996} and more
recently~\cite{Camporeale-Delzanno-Vergen-Moulton:2015}).
Although the choice of the Hermite basis is a crucial aspect of the
method (see also~\cite{Xiang-Wang:2013,Ma-Sun-Tang:2005}), appropriate
automatic algorithms, aimed to optimize the performances, are at the
moment not available.  
It should also be noted that Hermite-based method may suffer of
instability issues.
In alternative, the use of Legendre polynomials can offer stable
spectral approximations, with almost the same good conservation
properties of Hermite based discretizations
(see~\cite{Manzini:2016:JCP:journal,Manzini-Funaro-Delzanno:2017:journal}).

\smallskip

%%{\textcolor{red}{Aggiungero' una serie di references sui metodi spettrali dopo.
%%Il tutto comunque va riletto. Non mi piacciono alcuni punti.}}

\smallskip

The paper is organized as follows.
In Section~\ref{sec:continuous:Model}, we present the continuous
model, i.e., the 1D-1V Vlasov equation.
In Section~\ref{sec:phase-space discretization}, we introduce the
spectral approximation in the phase space.
In Section~\ref{sec:time:discretization}, we present the
semi-Lagrangian schemes based on a first-order accurate
approximation of the characteristic curves and we couple this with
a second-order BDF method.
In Section~\ref{sec:basis:choice} we examine Legendre and Hermite
discretizations with respect to the velocity $v$.
In Section~\ref{sec:conservation:properties}, we discuss the
conservation properties of the method in the discrete framework. 
%and define discrete quantities
%that represents the number of particles and the total energy.
%%
In Section~\ref{sec:numerical:experiments}, we assess the performance
of the method on a standard benchmark problem. In particular,
referring to the Hermite case, we make comparisons on the solution's
behavior depending on a parameter $\alpha$, which describes the decay
at infinity.
In Section~\ref{sec:conclusions}, we present our final remarks and
conclusions. We assume that the reader is confident with the main
results concerning spectral methods. Some passages in this paper are
given for granted, but they can easily recovered from standard texts,
such as, for instance,
\cite{Gottlieb-Orszag:1977,Canuto-Hussaini-Quarteroni-Zhang:1988,Funaro:1992,Boyd:1989,Shen-Tang-Wang:2011,Trefethen:2000,Guo:1998,Bernardi-Maday:1997}.

%{\textcolor{red}{Qui bisognerebbe aggiungere references generali sui metodi spettrali.
%Vedi Blocco 2}}  

%%%%%%%%%%%%%%%%%%%%%%%%%%%%%%%%%%%%%%%%%%%%%%%%%%%%%%%%%
\section{The 1D-1V  Vlasov equation}
\label{sec:continuous:Model}
%%%%%%%%%%%%%%%%%%%%%%%%%%%%%%%%%%%%%%%%%%%%%%%%%%%%%%%%%%

We consider the  Vlasov equation,
defined in a domain $\Omega =\Omega_x\times\Omega_v$, with $\Omega_x\subseteq\R$ and
$\Omega_v\subseteq\R$. The unknown $f=f(t,x,v)$, which denotes the probability of finding negative charged
particles at the location $x$ with velocity $v$,
 solves the problem:
\begin{align}
  &
  \frac{\partial\fs}{\partial t} + \vs\frac{\partial\fs}{\partial\xs} 
  -\Es(t,\xs)\ \frac{\partial\fs}{\partial\vs}=g (t,x,v), \qquad t\in (0,T],\,\,\xs\in\Omega_{x},\,\,
  \vs\in\Omega_{v}.
  \label{eq:1D1V:V}
\end{align}
Here the right-hand side $g$ is given. The function $f$ is required to satisfy the initial distribution:
\begin{align}
  & \fs(0,\xs,\vs)=\bar{\fs}(\xs,\vs),\quad\xs\in\Omega_{x},\,\,\vs\in\Omega_{v}.
  \label{eq:1D1V:Vci}
\end{align}
The electric field $\Es$ is coupled with $f$. Indeed, by setting:
\begin{equation}
  \label{eq:1D1V:EMaxwell}
  \frac{\partial\Es}{\partial\xs}(t,\xs) = 1-\rho(t,\xs),
  \quad t\in [0,T],\,\,\xs\in\Omega_{x},
\end{equation}
we have that the electron charge density $\rho(t,\xs)$ is defined
by:
\begin{equation}\label{eq:1D1V:charge-density}
  \rho(t,\xs)= \int_{\Omega_{v}}\fs(t,\xs,\vs)\dv.
  %% \quad t\in [0,T],\,\,\xs\in\Omega_{x}.
\end{equation}
In order to select a unique solution in~\eqref{eq:1D1V:EMaxwell},
we assume the following charge conservation property:
\begin{equation} 
  \label{eq:1D1V:intrhoconst}
  \int_{\Omega_{x}}\Es(t,\xs)dx = 0,
  \qquad \textrm{which~implies~that}\qquad 
  \int_{\Omega_{x}}\rho(t,\xs)dx = \abs{\Omega_x}, %\quad t\in [0,T],
\end{equation}
where $\abs{\Omega_x}$ is the size of $\Omega_x$.

\smallskip
Such a system of equations in the unknowns $f$ and $E$ 
is a simplification of the case where $\Omega_x$ and $\Omega_v$ are
two or three dimensional domains. In this extension  
the partial derivative in~\eqref{eq:1D1V:EMaxwell} is substituted by the divergence operator.
Thus, if the electric field takes the form of a gradient (i.e., $\vec E =\nabla \Phi$), the corresponding potential
must satisfy the Poisson equation $-\Delta \Phi =\rho$. The resulting set of equations is
the so called Vlasov-Poisson system. More advanced generalizations concern
with the coupling of the Vlasov equation with Maxwell's equations.

\smallskip
 
As far as boundary constraints are concerned, we will
assume periodic boundary conditions in the variable $x$  and
either periodic or homogeneous Dirichlet boundary conditions for the
variable $v$. In the case in which $\Omega_v=\R$, suitable exponential decay
conditions are assumed at infinity.

\smallskip

When $g=0$, from a straightforward calculation that takes into account boundary conditions, it follows that:
\begin{equation} 
  \label{eq:1D1V:massConserv}
  \frac{d}{dt}\int_{\Omega}\fs(t,\xs,\vs)\,dx\,dv = 0. %,\qquad t\in [0,T].
\end{equation}
This property is known as ``mass conservation''.
Moreover, after introducing the total energy of the system defined by:
\begin{equation}
  \label{eq:1D1V:totEnergy}
  \mathcal{E}(t) 
  = \frac12 \int_{\Omega}\fs(t,\xs,\vs)\,{|\vs|^{2}}\,dx\,dv\, 
  + \,\frac12\int_{\Omega_{x}}|\Es(t,\xs)|^{2 }\,dx,
  % ,\qquad t\in [0,T],
\end{equation}
where the first term represents the kinetic energy  
and the second one the potential energy, we recover
this other conservation property:
\begin{equation}
  \label{eq:1D1V:energyConserv}
  \frac{d}{dt}\,\mathcal{E} (t)=0.  %,\qquad t\in [0,T].
\end{equation}
If the electric field is smooth enough, for a sufficiently small
$\delta>0$, the local system of characteristics associated
with~\eqref{eq:1D1V:V} is given by the curves
$(\Xs(\tau),\Vs(\tau))$ solving: 
\begin{equation}
  \label{eq:1D1V:char}
  \frac{d\Xs}{d \tau} = -\Vs(\tau),\qquad
  \frac{d\Vs}{d \tau} =  \Es(\tau,\Xs(\tau)), \qquad \tau \in ]t-\delta ,t+\delta [,
\end{equation}
with the condition that $(\Xs(t),\Vs(t))=(\xs,\vs)$ when $\tau=t$.
Under suitable regularity assumptions, there exists a unique solution
of the Vlasov-Poisson problem \eqref{eq:1D1V:V}, \eqref{eq:1D1V:Vci},
\eqref{eq:1D1V:EMaxwell} and \eqref{eq:1D1V:charge-density},
(see, e.g., \cite{Glassey:1996}), which is formally expressed by
propagating the initial condition~\eqref{eq:1D1V:Vci} along the
characteristic curves described by~\eqref{eq:1D1V:char}.
Under hypotheses of smoothness, for every $t\in(0,T]$ we have that:
\begin{equation}
  \label{eq:1D1V:solChar}
  \fs(t,\xs,\vs) = \bar{\fs}(\Xs(t),\Vs(t)),
  %% \qquad t\in[0,T],\,\,\xs\in\Omega_{x},\,\,\vs\in\Omega_{v}.
\end{equation}
where we recall that $\bar f$ is the initial datum (see \eqref{eq:1D1V:Vci}).
By using the first-order approximation:
\begin{equation}
  \label{eq:1D1V:char1}
  \Xs(\tau) = \xs - \vs (\tau-t),\qquad
  \Vs(\tau) = \vs + \Es(t,\xs)(\tau-t) ,
\end{equation}
the Vlasov equation is satisfied up to an error that decays as
$\vert \tau-t \vert$, for $\tau$ tending to $t$.
%%

%%%%%%%%%%%%%%%%%%%%%%%%%%%%%%%%%%%%%%%%%%%%%%%%%%%%%%%%%%
\section{Phase-space discretization}
\label{sec:phase-space discretization}
%%%%%%%%%%%%%%%%%%%%%%%%%%%%%%%%%%%%%%%%%%%%%%%%%%%%%%%%%%
%
We set up a semi-Lagrangian spectral-type method to find numerical approximations
to the 1D-1V Vlasov-Poisson problem
given by equations \eqref{eq:1D1V:V}, \eqref{eq:1D1V:Vci},
\eqref{eq:1D1V:EMaxwell}  and \eqref{eq:1D1V:charge-density}.
The same algorithm was proposed in~\cite{Fatone-Funaro-Manzini:2018} within the framework of 
trigonometric functions, assuming a periodic distribution function $f$ in both $x$ and $v$ variables.
Instead, we will consider here a periodic function in $x$, and we will examine different sets
of basis functions in the variable $v$ and discuss pros and cons of the various approaches. 
The extension to higher-dimensional problems (the 3D-3V case, for instance) is
straightforward, though technically challenging in the implementation. 
The basic set up is discussed again in~\cite{Fatone-Funaro-Manzini:2018}.
\smallskip

Imposing periodic boundary conditions for the variable $x$ leads us
to consider the domain: $\Omega_{x} = [0,2\pi[$.
Given the positive integer $N$, we then consider the nodes:
\begin{equation}
  \xs_{i} = \frac{2\pi}{N}\,i ,\,\,\quad i=0,1,\ldots, N-1. 
\label{eq:1D1V:collpoints}
\end{equation}
Similarly, regarding the direction $v$, we will assume to have a set
of nodes:  $\vs_{j}, j=0,1,\ldots, M-1$, where $M$ is a given positive integer.
These nodes are explicitly defined later on.
Hereafter,  we  use the indices
$i$ and $n$ running from $0$ to $N-1$ to label the grid points along
the $x$-direction, and $j$ and $m$ running from $0$ to $M-1$ to label
the grid points along the $v$-direction.

\smallskip
We introduce the  Lagrangian
basis functions for the $\xs$ and $\vs$ variables with respect to
the nodes, that is:

\begin{equation}
  \Bs_{i}^{(N)} (\xs_{n}) = \delta_{in} \qquad\textrm{and}\qquad\Bs_{j}^{(M)}(\vs_{m}) = \delta_{jm},
  %%\Bs_{i}^{(N)} (\xs_{n}) = \delta_{in},\quad i,n=0,1,\ldots,N-1,\qquad\Bs_{j}^{(M)}(\vs_{m}) = \delta_{jm},\quad j,m=0,1,\ldots, M-1,
  \label{eq:1D1V:Kron}
\end{equation}
where $\delta_{ij}$ is the usual Kronecker symbol.
In the periodic case, it is known that:
\begin{align}
  \Bs_{i}^{(N)}(\xs) &= \frac{1}{N} \sin\left(\frac{N(\xs-\xs_{i})}{2}\right)\,\cot\left(\frac{\xs-\xs_{i}}{2}\right).\label{eq:1D1V:Basis1}
\end{align}
As far as the basis functions $\Bs_{j}^{(M)}$ are concerned, we will specify at due time
the various choices, according to the system of nodes adopted.
Furthermore, we define the discrete spaces:
\begin{align}
  \Xv_N    &= \textrm{span}\Big\{ \Bs_{i}^{(N)} \Big\}_{i=0,1,\ldots, N-1}
  \quad
  %%\nonumber\\[3mm]
  %%\qquad\xs\in\Omega_{x},\nonumber\\[3mm]
  \Yv_{N,M} = \textrm{span}\Big\{ \Bs_{i}^{(N)}\Bs_{j}^{(M)} \Big\}_{i=0,1,\ldots, N-1\atop j=0,1,\ldots, M-1}.
  %%\qquad\xs\in\Omega_{x},\,\vs\in\Omega_{v}.
  \label{eq:1D1V:spaces}
\end{align}
In this way, any function $\fs_{N,M}$ that belongs to $\Yv_{N,M}$ can
be represented as:
\begin{equation}
  \fs_{N,M}(\xs,\vs) = \sum_{i=0}^{N-1}\,\sum_{j=0}^{M-1} \cs_{ij}\,\Bs_{i}^{(N)}(\xs)\,\Bs_{j}^{(M)}(\vs),
  %%\qquad\xs\in\Omega_{x},\,\,\vs\in\Omega_{v},
  \label{eq:1D1V:fbasis}
\end{equation}
where the coefficients of such an expansion are given by:
\begin{equation}
  \cs_{ij} = \fs_{N,M} (\xs_{i},\vs_{j}).
  %%\quad i=0,1,\ldots,N-1,\,\,j=0,1,\ldots, M-1.
  \label{eq:1D1V:coeff}
\end{equation}

As usual, the matrix  $\ds_{ni}^{(N,s)}$ will denote the $s$-th derivative
of $\Bs_{i}^{(N)}$ evaluated at point $\xs_n$. Explicitly, we have:
\begin{equation}
  \ds_{ni}^{(N,s)}=\frac{d^{s}\Bs_{i}^{(N)}}{d \xs^{s}}(\xs_{n}) .
  \label{eq:1D1V:Der1v}
\end{equation}
In the same way, with respect to the variable $v$,  one defines $\ds_{mj}^{(M,s)}$.
As a special case we set: $\ds_{ni}^{(N,0)}=\delta_{ni}$, \
$\ds_{mj}^{(M,0)}=\delta_{mj}$.
%%
%Moreover, it is easy to prove that there exists a constant $C$,
%independent of $N$, such that:
%\begin{equation}
%  \abs{ \ds_{ni}^{(N,1)} } \leq CN.
%  \label{eq:1D1V:stimad}
%\end{equation}
%This estimate will be useful in the next section for studying the
%stability conditions in the time-marching schemes.

\smallskip 
Furthermore, in the trigonometric case, we remind the
following Gaussian quadrature formula (see, e.g.,
\cite{Canuto-Hussaini-Quarteroni-Zhang:1988} and
\cite{Sheng-Tao-Wang:2011}):
\begin{equation}\label{eq:1D1V:fdq}
  \displaystyle
  \frac{1}{2\pi}\,\int_{0}^{2\pi}\phi(\xs)\,d\xs\ \simeq \ \frac1N\sum_{i=0}^{N-1}\phi(\xs_{i}),
\end{equation}
which is exact for every $\phi$ belonging to the space:
$$\textrm{span}\Big\{1, \ \big\{\sin nx,\cos nx\big\}_{n=1,\ldots,N-1}, \ \sin Nx\Big\}.$$
%%
%% {\textcolor{red}{anche Shen, Tao, Wang, sempre section 2.1.2}.

%It is clear that, with little modifications, we can handle Lagrangian
%basis of nonperiodic type. 
%%
%Among these, the most representative ones are constructed on Legendre
%or Chebyshev algebraic polynomials, or Hermite functions (i.e.,
%Hermite polynomials multiplied by a Gaussian function).
%%
%In some preliminary tests, we observed that each one of these cases
%presents peculiar behavior in applications.
%%
\smallskip

Now, let us assume that the one-dimensional function $\Es_N\in\Xv_N$
is known. 
% is given.
%%
Given $\Delta t>0$, by taking $\tau=t -\Delta t$ in formula \eqref{eq:1D1V:char1}, we define the new
set of points $\{(\Xnewnm,\Vnewnm)\}_{n,m}$ where:
\begin{align}
  \Xnewnm &= \Xoldnm-\incrXnm      ,    \label{eq:1D1V:puntimossi}\\[2mm]
  \Vnewnm &= \Voldnm+E_N(x_n)\Delta t.  \label{eq:1D1V:1charnm}
\end{align}
To evaluate a function $\fs_{N,M}\in\Yv_{N,M}$ at the new points
$(\Xnewnm, \Vnewnm)$ through the coefficients in
\eqref{eq:1D1V:coeff}, we use a Taylor expansion in time.
%%
% If a function $\fs_{N,M}\in\Yv_{N,M}$ is known through its
% coefficients \eqref{eq:1D1V:coeff}, we would like to evaluate it at
% the new points $(\Xnewnm, \Vnewnm)$.
% %%
% To this end we use a Taylor expansions.
%%
%%We recall that, for a sufficiently smooth function $\Psi$, one has:
%%
For a generic smooth function $\Psi$, we have that
\begin{align}\label{eq:1D1V:Taylor}
  &\Psi (\xs-\vs\Delta t,\vs + \Es_N(\xs) \Delta t) 
  = \Psi(\xs,\vs) - \vs\Delta t\,\frac{\partial\Psi}{\partial\xs}(\xs,\vs) 
  + \Es_N(x) \Delta t\,\frac{\partial\Psi}{\partial\vs}(\xs,\vs)  + \ldots
\end{align}
where we omitted the terms in $\Delta t$ of order higher than one.
%% 
%% where  the usual little-$o$ notation  is used.\\
%% 
In particular, for  $\Psi(\xs,\vs)=\Bs_{i}^{(N)}(\xs)\,\Bs_{j}^{(M)}(\vs)$, one gets:
\begin{align}\label{eq:1D1V:eqclou2}
  &
  \Bs_{i}^{(N)}(\Xnewnm)\,\Bs_{j}^{(M)}(\Vnewnm) =
  \Bs_{i}^{(N)}(\Xoldnm)\,\Bs_{j}^{(M)}(\Voldnm) - 
  \incrXnm\, 
  \left[
    \frac{d \Bs_{i}^{(N)}}{dx}(\Xoldnm)
  \right]\,\Bs_{j}^{(M)}(\Voldnm)\nonumber\\[3mm]
  &\hskip1truecm		
  + \incrVnmb\,\Bs_{i}^{(N)}(\Xoldnm)\,\left[\frac{d \Bs_{j}^{(M)}}{d\vs}(\Voldnm)\right] + \ldots
\end{align}
Adopting the notation in \eqref{eq:1D1V:Kron} and \eqref{eq:1D1V:Der1v} the above relation
is rewritten as:
\begin{align}\label{eq:1D1V:eqclou3}
  & \Bs_{i}^{(N)}(\Xnewnm)\,\Bs_{j}^{(M)}(\Vnewnm) =
  \delta_{in}\,\delta_{jm}  - 
  \incrXnm\,\delta_{jm}\,\ds_{ni}^{(N,1)} +	
  \incrVnmb\,\delta_{in}\,\ds_{mj}^{(M,1)} + \ldots
\end{align}
Finally, by substituting \eqref{eq:1D1V:eqclou3} in \eqref{eq:1D1V:fbasis}, we
obtain the approximation:
\begin{align}\label{eq:1D1V:eqclouDef}
  \fs_{N,M}(\Xnewnm,\Vnewnm) 
  & =\sum_{i=0}^{N-1}\,\sum_{j=0}^{M-1}\cs_{ij}\,\Bs_{i}^{(N)}(\Xnewnm)\,\Bs_{j}^{(M)}(\Vnewnm)\nonumber\\[3mm]
  & \approx \sum_{i=0}^{N-1}\,\sum_{j=0}^{M-1}\cs_{ij} 
  \Big(
  \delta_{in}\,\delta_{jm} -  
  \incrXnm \,\delta_{jm}\,\ds_{ni}^{(N,1)} +	
  \incrVnmb\,\delta_{in}\,\ds_{mj}^{(M,1)}\Big)\nonumber\\
  & = \cs_{nm} +\Delta t\left[ -\vs_m\sum_{i=0}^{N-1}\ds_{ni}^{(N,1)}\cs_{im} 
    + \Es_N(x_n)\sum_{j=0}^{M-1}\ds_{mj}^{(M,1)}\cs_{nj}\right].\nonumber\\
\end{align}
We will see how to use the above relation in the following section.
%% 

%%%%%%%%%%%%%%%%%%%%%%%%%%%%%%%%%%%%%%%%%%%%%%%%%%%%%%%%%%
\section{Full discretization of the Vlasov equation}
\label{sec:time:discretization}
%%%%%%%%%%%%%%%%%%%%%%%%%%%%%%%%%%%%%%%%%%%%%%%%%%%%%%%%%%
Given the time instants $t^{k}=k \Delta t= k\,{T}\slash{K}$ for any
integer $k=0,1,\ldots,K$, we consider an approximation of
the unknowns $f$ and $E$  of problem
\eqref{eq:1D1V:V}, \eqref{eq:1D1V:Vci}, \eqref{eq:1D1V:EMaxwell},
\eqref{eq:1D1V:charge-density}. To this end, we set:
\begin{align}\label{eq:1D1V:fapproxt}
  \left(
    \fs_{N,M}^{(k)}(\xs,\vs),\,\Es^{(k)}_N(\xs) 
  \right) \simeq 
  \left(
    \fs(t^{k},\xs,\vs),\,\Es(t^{k},\xs) 
  \right) ,
  %,\nonumber \\[3mm]& \hskip5truecm  k=0,1,\ldots,K,\,\,
  \qquad\xs\in\Omega_{x},\,\,\vs\in\Omega_{v},
\end{align}
where the function $\fs_{N,M}^{(k)}$ belongs to $\Yv_{N,M}$ and the
function $\Es_N^{(k)}$ belongs to $\Xv_N$.
By taking into account \eqref{eq:1D1V:charge-density}, we define:
\begin{equation}\label{eq:1D1V:charge-densityapproxt}
  \displaystyle
  \rho^{(k)}_N(\xs)  =
  \int_{\Omega_{v}}\fs^{(k)}_{N,M}(\xs,\vs)\,\dv
  \ \simeq \ \rho(t^k,\xs).
  %%\qquad k=0,1,\ldots,K,\,\,\xs\in\Omega_{x}.
\end{equation}
Hence, at any time step $k$, 
%$k=0,1,\ldots, K$, 
we express $\fs^{(k)}_{N,M}$ in the following way:
\begin{align}\label{eq:1D1V:fbasisk}
  \fs^{(k)}_{N,M} (\xs,\vs) 
  = \sum_{i=0}^{N-1}\,\sum_{j=0}^{M-1} \cs^{(k)}_{ij}\,\Bs_{i}^{(N)}(\xs)\,\Bs_{j}^{(M)}(\vs) , 
  %\qquad k=0,1,\ldots, K,\,\,
  %\qquad\xs\in\Omega_{x},\,\,\vs\in\Omega_{v},
\end{align}
where
\begin{equation}\label{eq:1D1V:coeffK}
  \cs^{(k)}_{ij} = \fs^{(k)}_{N,M} (\xs_{i},\vs_{j}).
  %% , \qquad k=0,1,\ldots,K,\,\, 
 % \qquad i=0,1,\ldots, N-1,\,\,j=0,1,\ldots, M-1.
\end{equation}
At time $t=0$, we use the initial condition for $f$
(see \eqref{eq:1D1V:Vci}) by setting:
\begin{equation}\label{eq:1D1V:coeff0}
  \cs^{(0)}_{ij} = \fs(0,\xs_{i},\vs_{j})=\bar{\fs}(\xs_{i},\vs_{j}).
  %% , \qquad  i=0,1,\ldots, N-1,\,\,j=0,1,\ldots, M-1.
\end{equation}

We first examine the case $g=0$ and suppose that $\Es_N^{(k)}$ is given at step $k$. 
To this purpose we use standard procedures in spectral Fourier methods. According to  \cite{Fatone-Funaro-Manzini:2018}, we write:
\begin{equation}
  \Es^{(k)}_N(\xs) = 
  -\sum_{n=1}^{N/2} \frac{1}{n}
  \left[ 
    \hat{\as}_{n}^{(k)}\,\sin(n\xs) - \hat{\bs}_{n}^{(k)}\,\cos(n\xs)
  \right],
  %% \quad x\in\Omega_{x} = [0,2\pi[,
  \label{eq:1D1V:FourierSeriesE}
\end{equation}
where the discrete Fourier coefficients $\hat{\as}_{n}^{(k)}$ and
$\hat{\bs}_{n}^{(k)}$, $n=1,2,\ldots,N/2$,  are suitably related to those of $  \rho^{(k)}_N$. In particular, the behavior in time of the coefficient $\hat{\as}_{1}^{(k)}$ will be taken into consideration in the experiment's section.

By taking $\tau=t-\Delta t$ in \eqref{eq:1D1V:char1}, we define:
\begin{align}
  \Xnewnm &= \Xoldnm - \incrXnm,  \nonumber\\[2mm]
  \Vnewnm &= \Voldnm + \Es_N^{(k)}(\xs_n)\Delta t.
  %% ,  \nonumber\\[2mm]
  %% &\hskip2truecm \quad k=0,1,\ldots,K,\,\,n=0,1,\ldots, N-1,\,\,m=0,1,\ldots, M-1. 
  \label{eq:1D1V:2charnm}
\end{align}

Since, for $g=0$, $\fs$ is expected to
remain constant along the characteristics, the most straightforward discretization method
is then obtained by advancing the coefficients of $f_{N,M}$ as
follows:
\begin{align}
  \cs^{(k+1)}_{nm} 
  = \fs^{(k)}_{N,M}(\Xnewnm,\Vnewnm) 
  = \sum_{i=0}^{N-1}\,\sum_{j=0}^{M-1}\cs^{(k)}_{ij}\,\Bs_{i}^{(N)}(\Xnewnm)\,\Bs_{j}^{(M)}(\Vnewnm),
  \label{eq:1D1V:favanz}
\end{align}
where we used ~\eqref{eq:1D1V:fbasisk}.
This states that %, for any time-step $k$, $k=0,1,\ldots, K$,
the value of $\fs_{N,M}^{(k+1)}$, at the grid points and time step
$(k+1)\Delta t$, is assumed to be equal to the previous value at time
$k\Delta t$, recovered by going backwards along the characteristics.
Technically, in \eqref{eq:1D1V:2charnm} we should use
$\Es_N^{(k+1)}(\xs_n)$ instead of $\Es_N^{(k)}(\xs_n)$, thus arriving
at an implicit method.
However, the distance between these two quantities is of the order of
$\Delta t$, so that the replacement has no practical effects on the
accuracy of  first-order methods. As a matter of fact, by computing the direction of the characteristic lines according to
\eqref{eq:1D1V:2charnm}, the scheme we are going to discuss turns out to be only first-order
accurate in $\Delta t$.
For higher order schemes, things must be treated more carefully.

\smallskip
Between each step $k$ and the successive one, we need to update the
electric field.
This can be done as follows.
Let $t^{k}$ be fixed.
%%Let $t^{k}$, $k=0,1,\ldots, K$, be fixed.
%% 
Using the Gaussian quadrature formula \eqref{eq:1D1V:fdq} in
\eqref{eq:1D1V:charge-densityapproxt}, and, then, applying~\eqref{eq:1D1V:coeffK} we
write:
\begin{equation}\label{eq:1D1V:charge-densityapproxt2}
  \displaystyle
  \rho^{(k)}_N(\xs_{i})
  = \frac{2\pi}{M}\,\sum_{j=0}^{M-1}\,\fs^{(k)}_{N,M}(\xs_{i},\vs_{j})
  = \frac{2\pi}{M}\,\sum_{j=0}^{M-1}\,\cs_{ij}^{(k)},
  %\qquad i=0,1,\ldots, N-1.
\end{equation}
where $\xs_{i}$, $i=0,\ldots,N-1$, are the nodes introduced in~\eqref{eq:1D1V:collpoints}.
At this point, in order to compute the new point-values  $\Es_N^{(k+1)}(\xs_n)$ of 
the electric field, we need to integrate  $\rho^{(k)}_N$.
To this purpose we use standard procedures in spectral Fourier methods.
 Note that the process is not
so simple in the higher dimensional case, where the operator to be inverted is the divergence.
For this reason, it turns out to be more convenient to solve a Poisson problem
for the potential of the electric field. This is also a very consolidated procedure
within the framework of spectral techniques.
%\CITE{}.
%%
\smallskip

By substituting~\eqref{eq:1D1V:eqclou3} in~\eqref{eq:1D1V:favanz}, we end up with the scheme:
\begin{align}
  \cs^{(k+1)}_{nm} 
  = \cs_{nm}^{(k)} + \Delta t\,\Phi_{nm}^{(k)},
  \label{eq:1D1V:favanza}
\end{align}
% \begin{align}
%   \cs^{(k+1)}_{nm} 
%   &= \cs_{nm}^{(k)} + \Delta t\,\Phi_{nm}^{(k)} \ ,\nonumber\\[3mm]
%   &\hskip1truecm	
%   \quad k=0,1,\ldots,K-1,\,\,n=0,1,\ldots, N-1,\,\,m=0,1,\ldots, M-1,
%   \label{eq:1D1V:favanza}
% \end{align}
where
\begin{align}\label{eq:1D1V:favanzadef}
  \Phi_{nm}^{(k)} 
  =
  - \vs_m\sum_{i=0}^{N-1}\ds_{ni}^{(N,1)}\cs_{im}^{(k)} 
  + \Es_N^{(k)}(\xs_n)\sum_{j=0}^{M-1}\ds_{mj}^{(M,1)}\cs_{nj}^{(k)}.
\end{align}
%%
% \begin{align}\label{eq:1D1V:favanzadef}
%   \Phi_{nm}^{(k)} 
%   &=
%   - \vs_m\sum_{i=0}^{N-1}\ds_{ni}^{(N,1)}\cs_{im}^{(k)} 
%   + \Es_N^{(k)}(\xs_n)\sum_{j=0}^{M-1}\ds_{mj}^{(M,1)}\cs_{nj}^{(k)} \nonumber\\
%   & \hskip1truecm	
%   \quad k=0,1,\ldots, K-1,\,\,n=0,1,\ldots, N-1,\,\,m=0,1,\ldots, M-1.
% \end{align}

%% In view of solving the non-homogeneous 1D-1V Vlasov-Poisson equation:

%%
%To solve the non-homogeneous Vlasov equation:
%%% 
%\begin{equation}\label{eq:1D1V:vlasovno}
%  \displaystyle
%  \frac{\partial\fs}{\partial t} + \vs\ \frac{\partial\fs}{\partial\xs} 
%  - \Es(t,\xs)\ \frac{\partial\fs}{\partial\vs} = \gs, 
%  %% \quad t\in (0,T],\,\,\xs \in \Omega_{x},\,\, v\in \Omega_{v},
%\end{equation}
%% 
%% where $\gs(t,\xs,\vs)$, $t\in(0,T]$, $\xs\in\Omega_{x}$,
%% $\vs\in\Omega_{v}$, is a given sufficiently regular right-hand side,

If we want to include a non-zero right-hand side term $\gs$, assumed to
be defined on $\Omega$ for every $t\in[0,T]$ and to be  sufficiently
regular, it is enough to modify \eqref{eq:1D1V:favanza} as follows:
\begin{align}\label{eq:1D1V:favanzano}
  %% & 
  \cs^{(k+1)}_{nm} =
  \cs_{nm}^{(k)} + \Delta t\,\Phi_{nm}^{(k)} + \Delta t\,\gs(t^k,\xs_n,\vs_m),
%   ,\nonumber\\[3mm]
%   &\hskip1truecm	
%   \quad k=0,1,\ldots, K-1,\,\,n=0,1,\ldots, N-1,\,\,m=0,1,\ldots, M-1.
\end{align}
where $\Phi_{nm}^{(k)}$
%% , $k=0,1,\ldots, K-1$, $ n=0,1,\ldots, N-1$, $m=0,1,\ldots, M-1$,
is the same as in \eqref{eq:1D1V:favanzadef}.
The scheme obtained is basically a forward Euler iteration. 
\smallskip

As expected from an explicit method, the parameter $\Delta t$ must
satisfy a suitable CFL condition, which is obtained by
requiring that the point $(\tilde{\xs}_{nm}, \tilde{\vs}_{nm})$ falls
inside the box $]x_{n-1}, x_{n+1}[\times ]v_{m-1}, v_{m+1}[$.  
From \eqref{eq:1D1V:2charnm}, a sufficient restriction is given by:
\begin{equation}\label{eq:1D1V:CFL}
  \Delta t 
  \leq 2\pi \Big( N \max_m \vert v_m\vert + M^\sigma \max_n \vert E^{(k)}_N(x_n)\vert\Big)^{-1}.
\end{equation}
The exponent $\sigma >0$ depends from the technique we use to handle numerically
the variable $v$. 
\smallskip 
A straightforward way to increase the time accuracy is to use a
multistep discretization scheme.
To this end, we consider the second-order accurate two-step 
Backward Differentiation Formula (BDF).
With the notation in \eqref{eq:1D1V:favanz},
\eqref{eq:1D1V:favanzadef} and \eqref{eq:1D1V:favanzano},  we have:
\begin{align}\label{eq:1D1V:favanzabdf2}
  %& 
  \fs_{N,M}^{(k+1)}(\xs_n,\vs_m) 
  = {{\frac43}} \fs_{N,M}^{(k)}  (\tilde{\xs}_{nm}, \tilde{\vs}_{nm}) 
  - {{\frac13}} \fs_{N,M}^{(k-1)}(\tilde{\tilde{\xs}}_{nm}, \tilde{\tilde{\vs}}_{nm}) 
  + {{\frac23}} \Delta t\,\gs(t^{k+1},\xs_n,\vs_m), 
%   \nonumber\\[3mm]
%   & \hskip1truecm	
%   \quad k=1,2,\ldots, K-1,\,\,n=0,1,\ldots, N-1,\,\,m=0,1,\ldots, M-1,
\end{align}
where, based on \eqref{eq:1D1V:2charnm},
$(\tilde{\xs}_{nm},\tilde{\vs}_{nm})$ 
%% $n=0,1,\ldots,N-1$, $m=0,1,\ldots,M-1$,
%%
is the point obtained from $(\xs_n,\vs_m)$ going
back of one step $\Delta t$ along the characteristic lines.
Similarly, 
%% for $n=0,1,\ldots,N-1$, $m=0,1,\ldots,M-1$,
the point $(\tilde{\tilde{\xs}}_{nm},\tilde{\tilde{\vs}}_{nm})$ is
obtained by going two steps back along the characteristic lines
(i.e., by replacing $\Delta t$ with $2\Delta t$ in \eqref{eq:1D1V:2charnm}).
Note that if $\gs=0$, it turns out that $\fs_{N,M}$ is constant along
the characteristic lines. Note also that a backward step of the amount of $2\Delta t$
has effect of the CFL condition, so that the restriction \eqref{eq:1D1V:CFL}
on the parameter $\Delta t$ should be halved. Despite the fact that a BDF scheme is
commonly presented as an implicit technique, for our special equation ($f$ constant
along the characteristics) it will assume the form of an explicit method.

%% By approximating at the first order the above values as follows:
%%
\smallskip
From the practical viewpoint, we can make the estimates:
\begin{align}\label{eq:1D1V:appf}
  & \fs_{N,M}^{(k)}  (\tilde{\xs}_{nm},\tilde{\vs}_{nm}) \ \simeq \ \cs_{nm}^{(k)} + \Delta t\,\Phi_{nm}^{(k)}, \nonumber\\[3mm]
  & \fs_{N,M}^{(k-1)}(\tilde{\tilde{\xs}}_{nm},\tilde{\tilde{\vs}}_{nm})\ \simeq \ \cs_{nm}^{(k-1)} + 2\Delta t\,\Phi_{nm}^{(k-1)}.
  % \nonumber\\[3mm]
%   & \hskip1truecm	
%   \quad k=1,2,\ldots, K-1,\,\,n=0,1,\ldots, N-1,\,\,m=0,1,\ldots, M-1,
\end{align}
Therefore, in terms of the coefficients, we end up with the scheme:
\begin{align}\label{eq:1D1V:favanzabdf2co}
  & \cs_{nm}^{(k+1)} =
  {\frac43} \Big( \cs_{nm}^{(k)}   +  \Delta t\,\Phi_{nm}^{(k)}  \Big) -
  {\frac13} \Big( \cs_{nm}^{(k-1)} + 2\Delta t\,\Phi_{nm}^{(k-1)} \Big) + {\frac23}\Delta t\,\gs(t^{k+1},\xs_n,\vs_m) \nonumber\\
  &\hskip1truecm = {{\frac43}} \cs_{nm}^{(k)} - {{\frac13}} \cs_{nm}^{(k-1)} 
  + {{\frac23}}\Delta t \left[  
    -\vs_m\sum_{i=0}^{N-1} \ds_{ni}^{(N,1)} (2\cs_{im}^{(k)} - \cs_{im}^{(k-1)}) \right. \nonumber\\
  & \hskip1truecm\qquad\left. + \Es_N^{(k)}(x_n)\sum_{j=0}^{M-1} \ds_{mj}^{(M,1)} (2\cs_{nj}^{(k)} - \cs_{nj}^{(k-1)})\right]+
  {\frac23}\Delta t\,\gs(t^{k+1},\xs_n,\vs_m).
%   ,  \nonumber\\[3mm]
%   & \hskip1truecm	
%   \quad k=1,2,\ldots, K-1,\,\,n=0,1,\ldots, N-1,\,\,m=0,1,\ldots, M-1.
\end{align}
From the experiments shown in~\cite{Fatone-Funaro-Manzini:2018}, it
turns out that this method is actually second-order accurate in
$\Delta t$. Higher order schemes can be obtained with similar
principles.

%%%%%%%%%%%%%%%%%%%%%%%%%%%%%%%%%%%%%%%%%%%%%%%%%%%%%%%%%
\section{Choice of the basis for the variable $v$}
\label{sec:basis:choice}
%%%%%%%%%%%%%%%%%%%%%%%%%%%%%%%%%%%%%%%%%%%%%%%%%%%%%%%%%%

Here, we will continue to assume that our problem is periodic with
respect to the variable $x$ and examine alternative techniques to
handle the variable $v$.
If we impose periodic boundary conditions for $v$, we will set:
$\Omega_{v} = [0,2\pi[$.
For a positive integer $M$, the nodes are:
\begin{equation}
  \vs_{j} = \frac{2\pi}{M}\,j ,\,\,\quad j=0,1,\ldots, M-1. 
\label{eq:1D1V:collpointsv}
\end{equation}
The Lagrangian basis functions are similar to those in
\eqref{eq:1D1V:Basis1}, that is:
\begin{align}
  \Bs_{j}^{(M)}(\vs) &= \frac{1}{M} \sin\left(\frac{M(\vs-\vs_{i})}{2}\right)\,\cot\left(\frac{\vs-\vs_{j}}{2}\right).\label{eq:1D1V:Basis1v}
\end{align}

\smallskip
Instead, if we prefer to work with algebraic polynomials, 
two choices are particularly suggested.
The first one is related to Legendre polynomials, the other one to
Hermite polynomials. In the first case, we set $\Omega_{v} =
]-1,1[$. The nodes $\vs_{j}$, $j=0,1,\ldots, M-1$ are the zeros of
$P^\prime_{M+1}$, where $P_{M+1}$ denotes the Legendre polynomial of
degree $M+1$.  In addition, we require that $ \Bs_{j}^{(M)}(\pm 1)=0$.
As a consequence, the Lagrangian basis functions take the form:
\begin{align}
  \Bs_{j}^{(M)}(\vs) &= \frac{(v^2-1)P^\prime_{M+1}(v)}{(M+1)(M+2)(v-v_j)P_{M+1}(v_j)}.
  \label{eq:1D1V:Basis1vleg}
\end{align}

For a fixed integer $s\geq 0$, explicit expressions of  $\ds_{mj}^{(M,s)}$
are well-known.
%\CITE{}
%% 
Moreover, we recall the Gaussian quadrature formula:
\begin{equation}\label{eq:1D1V:fdqleg}
  \displaystyle
  \int_{-1}^{1}\phi(\vs)\omega (v)\,d\vs \ \simeq \ \sum_{j=0}^{M-1}\phi(\vs_{j})w_j,
\end{equation}
which is exact for every $\phi$ belonging to the space of polynomials of
degree less or equal to $2M+1$ satisfying the condition $\phi (\pm 1)=0$. The weight function $\omega (v)=1$ is constant and the positive weights $w_j$, $j=0,\ldots, M-1$
are also well-known.
%\CITE{}
%%
This means that \eqref{eq:1D1V:charge-densityapproxt2} must be replaced by:
\begin{equation}\label{eq:1D1V:charge-densityapproxt2leg}
  \displaystyle
  \rho^{(k)}_N(\xs_{i})
  = \sum_{j=0}^{M-1}\,\fs^{(k)}_{N,M}(\xs_{i},\vs_{j}) w_j
  = \sum_{j=0}^{M-1}\,\cs_{ij}^{(k)} w_j . 
  % \qquad i=0,1,\ldots, N-1.
\end{equation}
where we assume that, in the variable $v$, $\fs^{(k)}_{N,M}$ is an algebraic 
polynomial of degree less or equal to $M+1$, such that $\fs^{(k)}_{N,M}(x_i,\pm 1)=0$.
\smallskip

In the Hermite case we have  $\Omega_v =]-\infty, +\infty[$, and the nodes
$\vs_{j}$, $j=0,\ldots, M-1$ are the zeros of $H_M$, which is the
Hermite polynomial of degree $M$.
The corresponding Lagrangian basis is:
\begin{align}
  \Bs_{j}^{(M)}(\vs) &= \frac{H_M(v)}{(v-v_j)H^\prime_M(v_j)} .\label{eq:1D1V:Basis1vher}
\end{align}

Also in this case, given the integer $s\geq 0$, explicit expressions
of $\ds_{mj}^{(M,s)}$ are well-known.
%\CITE{}
%%
The Gaussian quadrature formula is the same as in \eqref{eq:1D1V:fdq},
but now we have $\omega (v)=\exp (-v^2)$.
This is exact for every $\phi$ belonging to the space of polynomials
of degree less or equal to $2M+1$.
The positive weights $w_j$, $j=0,\ldots, M-1$ can be easily found in the literature.

%found for
%instance in~\cite{Canuto-Hussaini-Quarteroni-Zhang:1988,Guo:1998}.%%\CITE{}.

\smallskip 
In this context, it is advisable to make the change of variable
$f(t,x,v)=p(t,x,v)\exp(-v^2)$ in the Vlasov equation, so obtaining:
\begin{align}
  &
  \frac{\partial p}{\partial t} + \vs\frac{\partial p}{\partial\xs} 
  -\Es(t,\xs)\left[ \frac{\partial p}{\partial\vs} -2vp\right] =g (t,x,v)e^{v^2},\ 
  t\in (0,T],\,\,\xs\in\Omega_{x},\,\,
  \vs\in\Omega_{v}.
  \label{eq:1D1V:Vp}
\end{align}
At time step $k$, the function $p(t^k,x,v)$ is approximated by a
function $p^{(k)}_{N,M}(x,v)$ belonging to the finite dimensional
space $\Yv_{N,M}$.
Hence, with respect to the variable $v$, $ p^{(k)}_{N,M}$ is a
polynomial of degree at most $M-1$.
As a consequence, with abuse of notation, the coefficients $\cs^{(k)}_{ij}$ are modified as
follows:
\begin{equation}
  \cs^{(k)}_{ij} = p^{(k)}_{N,M} (\xs_{i},\vs_{j})= f^{(k)}_{N,M} (\xs_{i},\vs_{j})\, e^{v^2_j}.
  %% \quad i=0,1,\ldots,N-1,\,\,j=0,1,\ldots, M-1.
  \label{eq:1D1V:coeffher}
\end{equation}
% 
% \textcolor{red}{Forse qui conviene usare un simbolo diverso per i coefficienti di $p$ (finora $\cs^{(k)}_{ij}$ indicano i coefficienti di $f$)
where, now, relation \eqref{eq:1D1V:charge-densityapproxt2} is
replaced by:
\begin{equation}\label{eq:1D1V:charge-densityapproxt2her}
  \displaystyle
  \rho^{(k)}_N(\xs_{i})
  = \sum_{j=0}^{M-1}\,p^{(k)}_{N,M}(\xs_{i},\vs_{j}) w_j
  = \sum_{j=0}^{M-1}\,\cs_{ij}^{(k)} w_j , 
  % \qquad i=0,1,\ldots, N-1.
\end{equation}
while \eqref{eq:1D1V:favanzadef} is replaced by:
\begin{align}\label{eq:1D1V:favanzadefp}
  \Phi_{nm}^{(k)} 
  =
  - \vs_m\sum_{i=0}^{N-1}\ds_{ni}^{(N,1)}\cs_{im}^{(k)} 
  + \Es_N^{(k)}(\xs_n)\left[ \sum_{j=0}^{M-1}\ds_{mj}^{(M,1)}\cs_{nj}^{(k)}
-2v_m \cs_{nm}^{(k)}\right].
\end{align}

\smallskip
A straightforward generalization consists in introducing a parameter
$\alpha$ and assume that the weight function is $\omega (v)=\exp
(-\alpha^2 v^2)$.
This basically corresponds to a suitable stretching of the real axis
$\R$.
Thus, the approximation scheme can be easily adjusted by modifying
nodes and weights of the Gaussian formula, through a multiplication by
suitable constants. 
The difficulty in the implementation is practically the same, but, as
we shall see from the experiments, the results may be very sensitive
to the variation of $\alpha$.

\begin{figure}
  \centerline{
    \includegraphics[height=5.5cm]{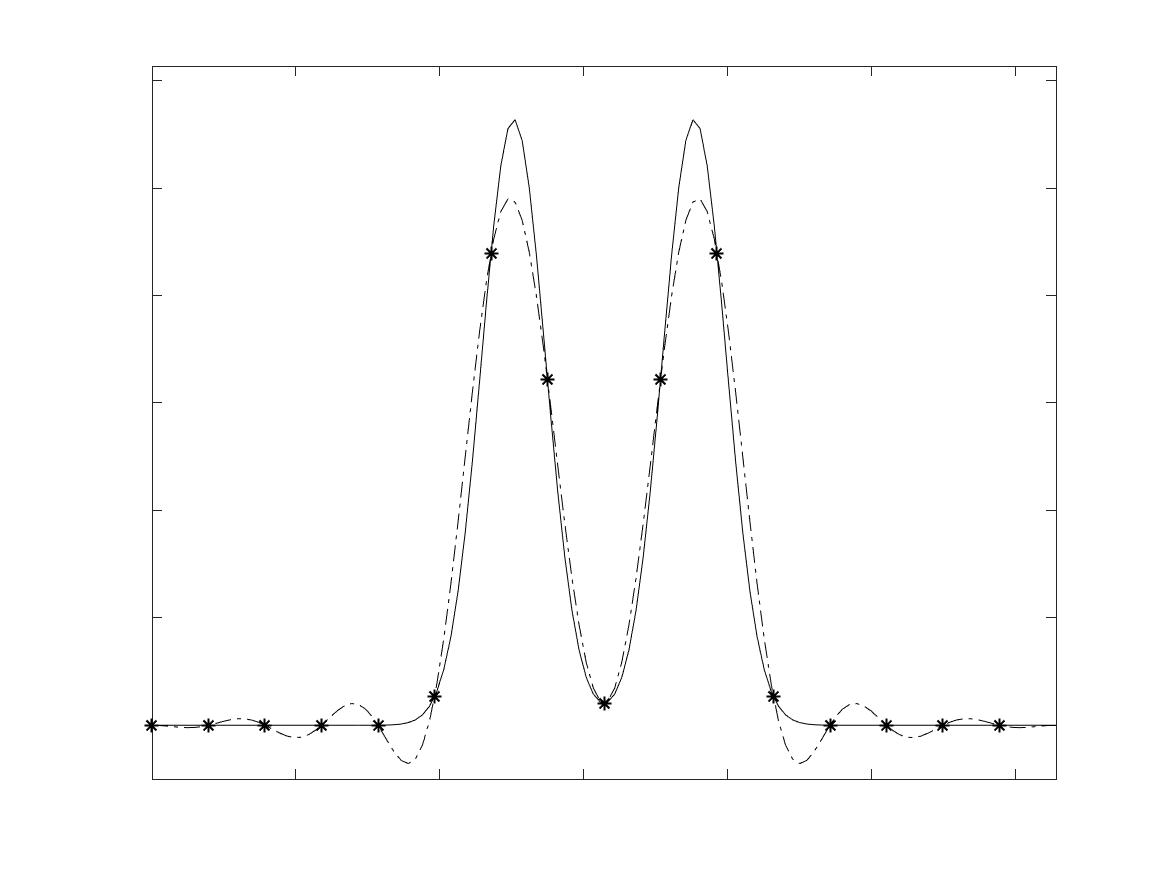}%{./figurePaper/GaussN2^4.jpg} 
    \hskip-0.6truecm
    \includegraphics[height=5.5cm]{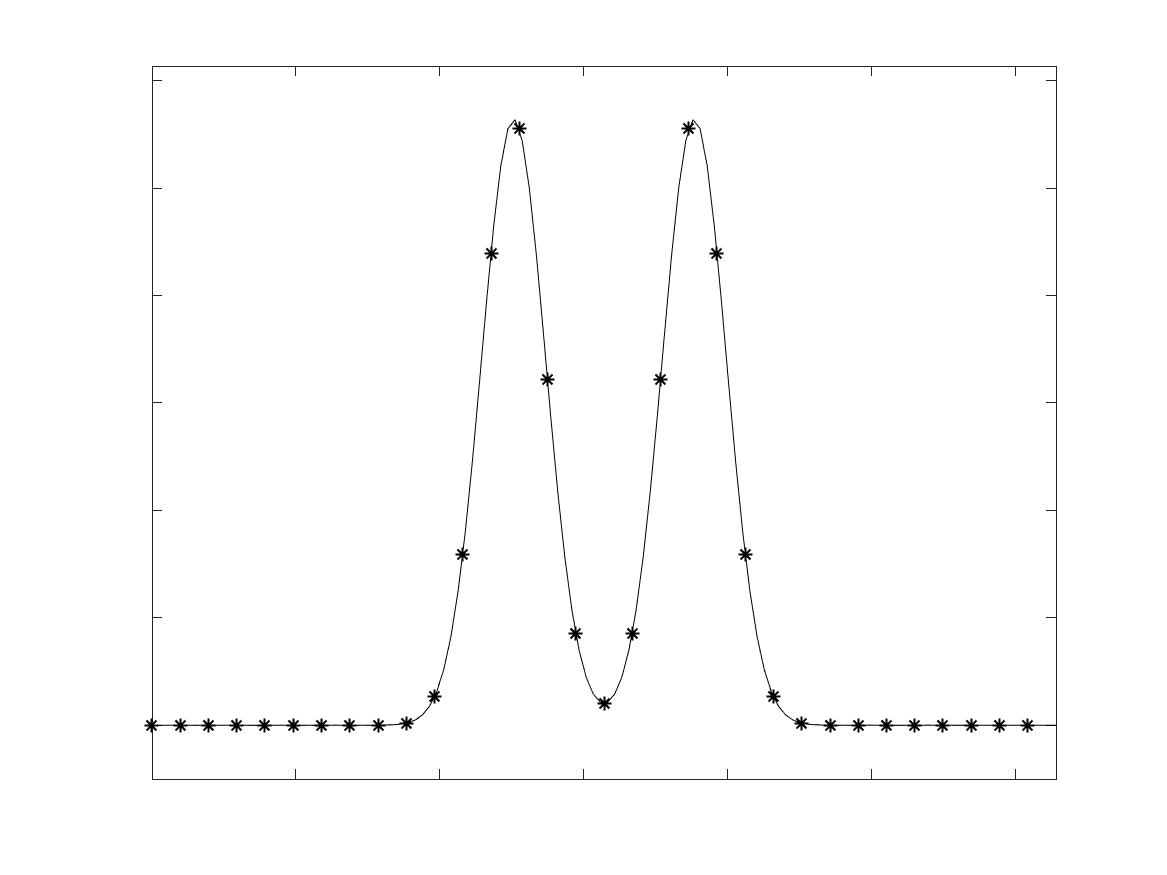}%{./figurePaper/GaussN2^5.jpg} 
  }
   \centerline{
    \includegraphics[height=5.5cm]{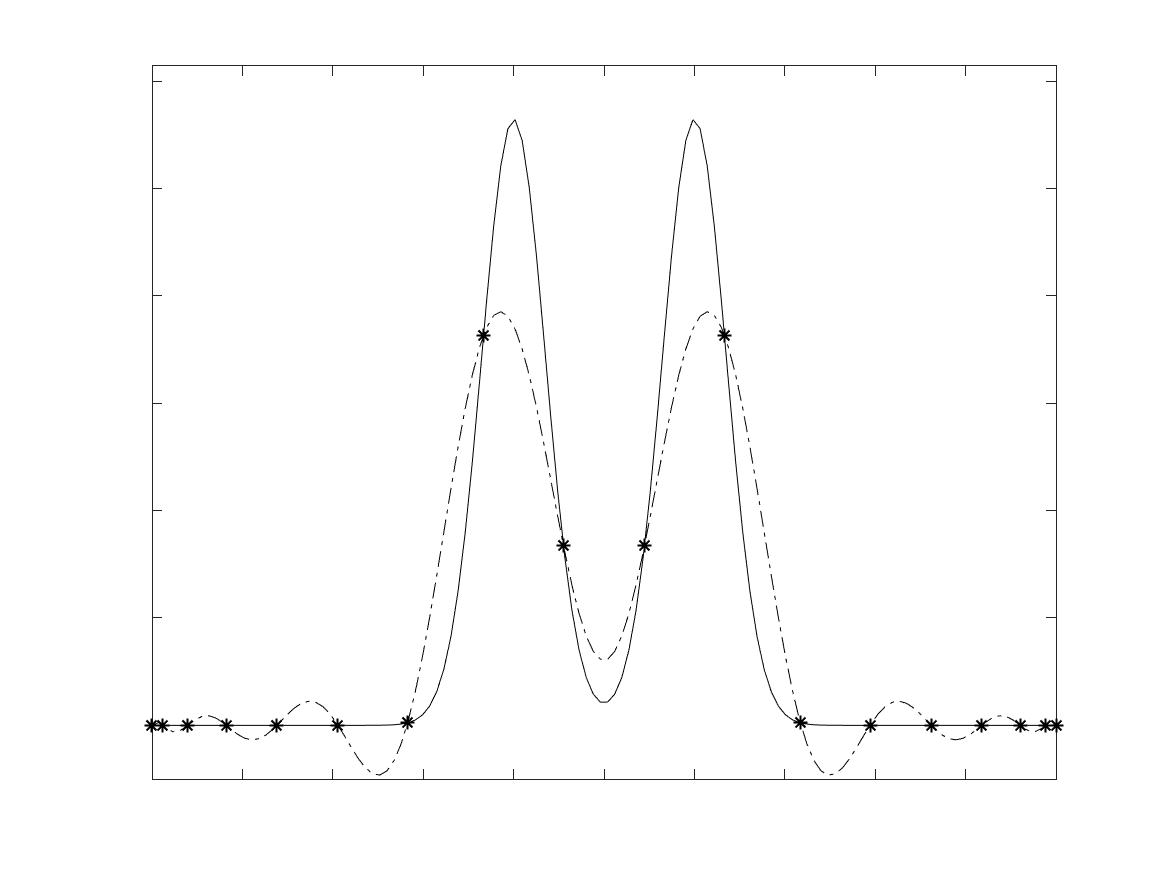}%{./figurePaper/LegendreN2^4.jpg} 
    \hskip-0.6truecm
     \includegraphics[height=5.5cm]{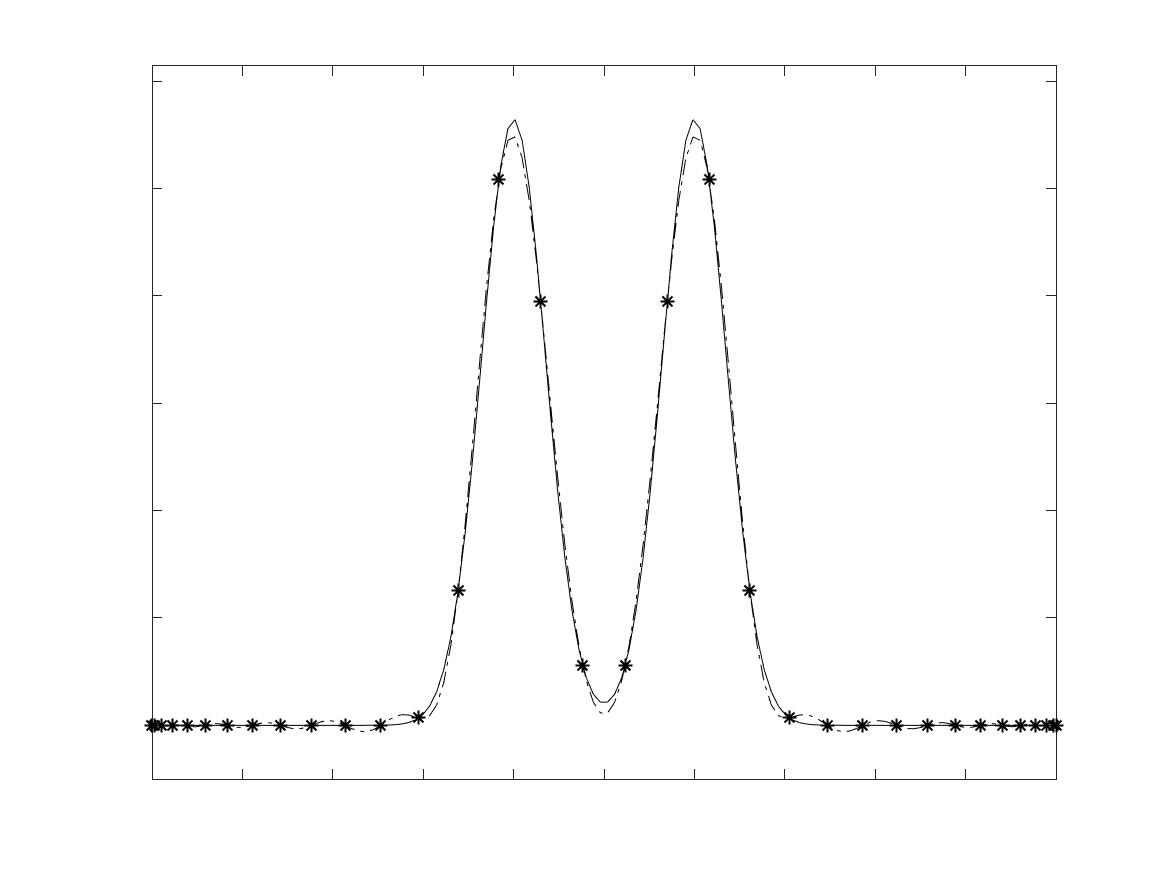}%{./figurePaper/LegendreN2^5.jpg} 
  }
  \centerline{
    \includegraphics[height=5.5cm]{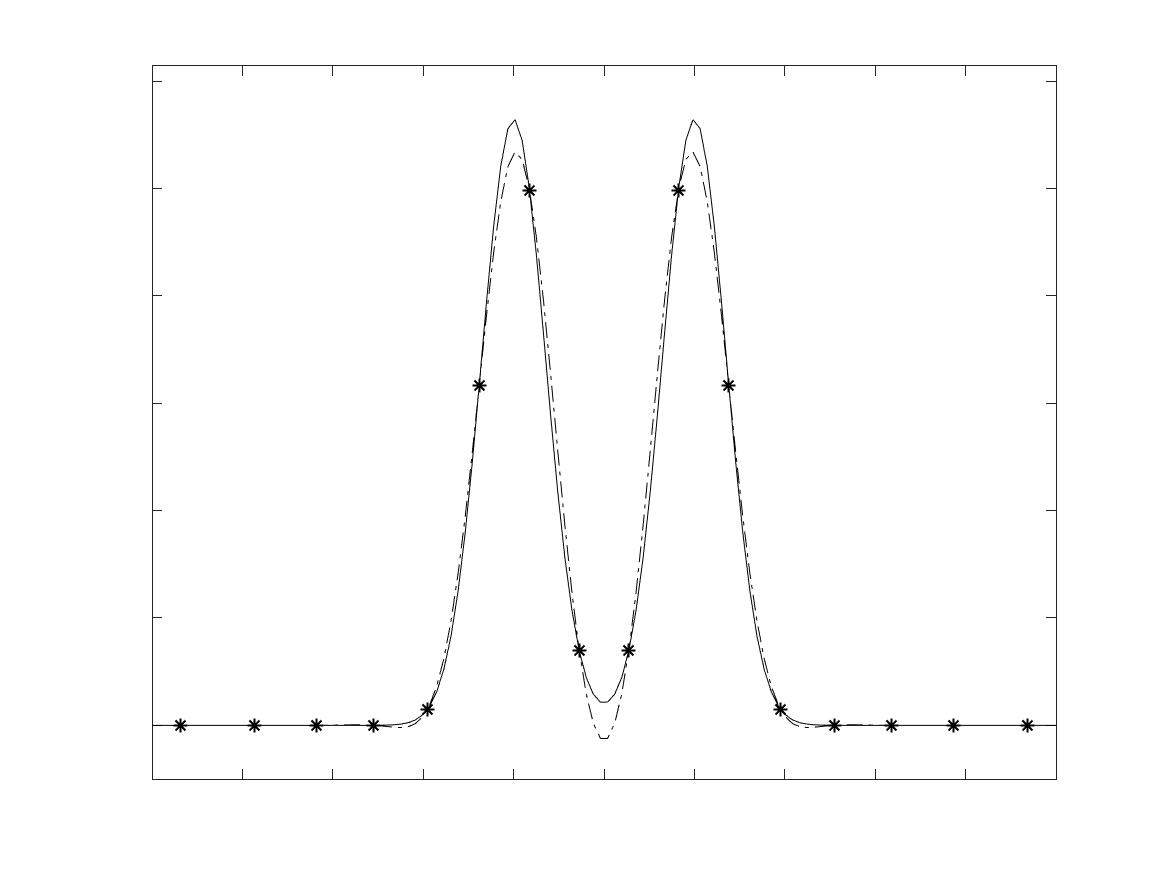}%{./figurePaper/HermiteN2^4.jpg} 
    \hskip-0.6truecm
    \includegraphics[height=5.5cm]{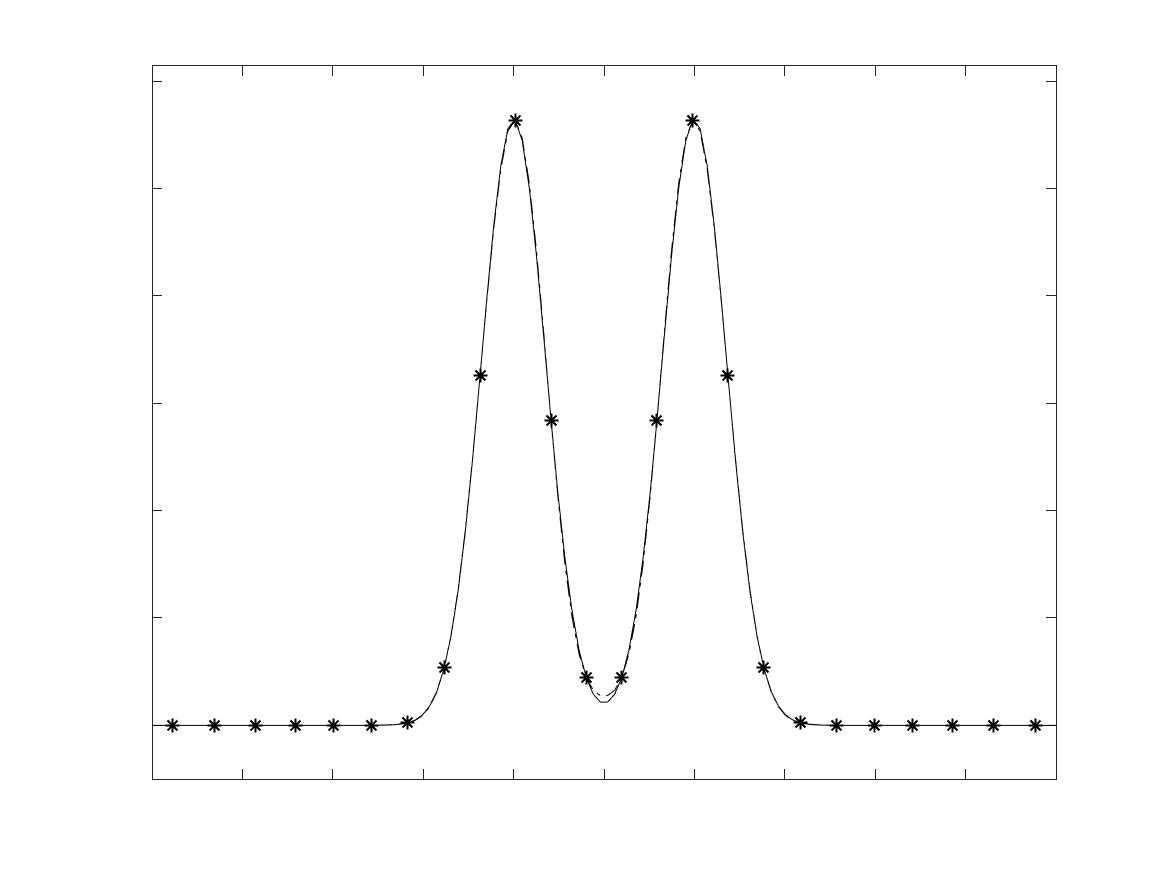}%{./figurePaper/HermiteN2^5.jpg} 
  }
  \caption{Interpolation with respect to $v\in [-5,5]$ of the function
    \eqref{eq:1D1V:twostreamfzero} with $M=2^{4}$ (left plots)
    or $M=2^{5}$ (right plots).  
   After a suitable mapping of the nodes into the interval $[-5,5]$,  the interpolations are relative 
   to the  Fourier case (top), the Legendre case (center) and the Hermite
    case (bottom).
  }
  \label{fig1approx}
\end{figure}

\smallskip
To see how the proposed ideas work in practice we turn our attention
to a special function to be interpolated, namely:
\begin{align} \label{eq:1D1V:twostreamfzero} 
   \xi(v)= \exp { \left(-\frac{\vs-\beta}{a\sqrt{2}}\right)^{2}}+
    \exp{ \left(-\frac{\vs+\beta}{a\sqrt{2}}\right)^{2}}
  \qquad \vs\in [-5,5],
\end{align}
with $a=1\slash{\sqrt 8}$, $\beta=1$.
The above is a multiple of the initial datum \eqref{eq:1D1V:twostreamf}, for $x=0$, used in the numerical
experiments of Section \ref{sec:numerical:experiments}.
We map the interval $[0,2\pi ]$ into the interval $[-5,5]$ in the case of Fourier type approximations,  
the interval $[-1, 1]$ into $[-5,5]$ in the case of Legendre type approximations. No mapping is made in the Hermite case.
In order to study the behavior with respect to the
variable $v$,  different sets of interpolation nodes are implemented.
Figure \ref{fig1approx} shows the results for $M=2^{4}$ (left plots)
and for $M=2^{5}$ (right plots).
On top, Fourier equispaced nodes \eqref{eq:1D1V:collpointsv} are
considered.
 Legendre and Hermite nodes are used for the plots in the middle and on
 bottom, respectively.
These last experiments correspond to the choice $\alpha =1$ for the
weight function $\exp (-\alpha^2 v^2)$. 

\smallskip

As far as the Legendre case is concerned, it has to be noted that
this is the one displaying the poorest interpolation properties. Such an observation turns out to be true in particular for $M=16$.
Despite the fact the Legendre polynomials are an excellent tool to approximate differential problems,
here the nodes accumulate at the extremes of the interval $]-1 ,1[$. From the physics viewpoint, this means that we are
concentrating too much attention to particles with high velocity $v$, that are usually less frequent.
As we shall see later, this property reflects negatively on the performances of Legendre approximations
when applied to the full Vlasov-Poisson  problem. Of course, one can always increase the number of nodes in the variable $v$.
However, at this stage, we are only observing that periodic Fourier or Hermite functions can represent a valid
alternative. We also point out that the Fourier approach opens the way to the use of fast algorithms,
such as the Discrete Fast Fourier Transform
(DFT)~\cite{Brigham:1988}.

\smallskip

\begin{figure}
  \centerline{
    \includegraphics[height=5.5cm]{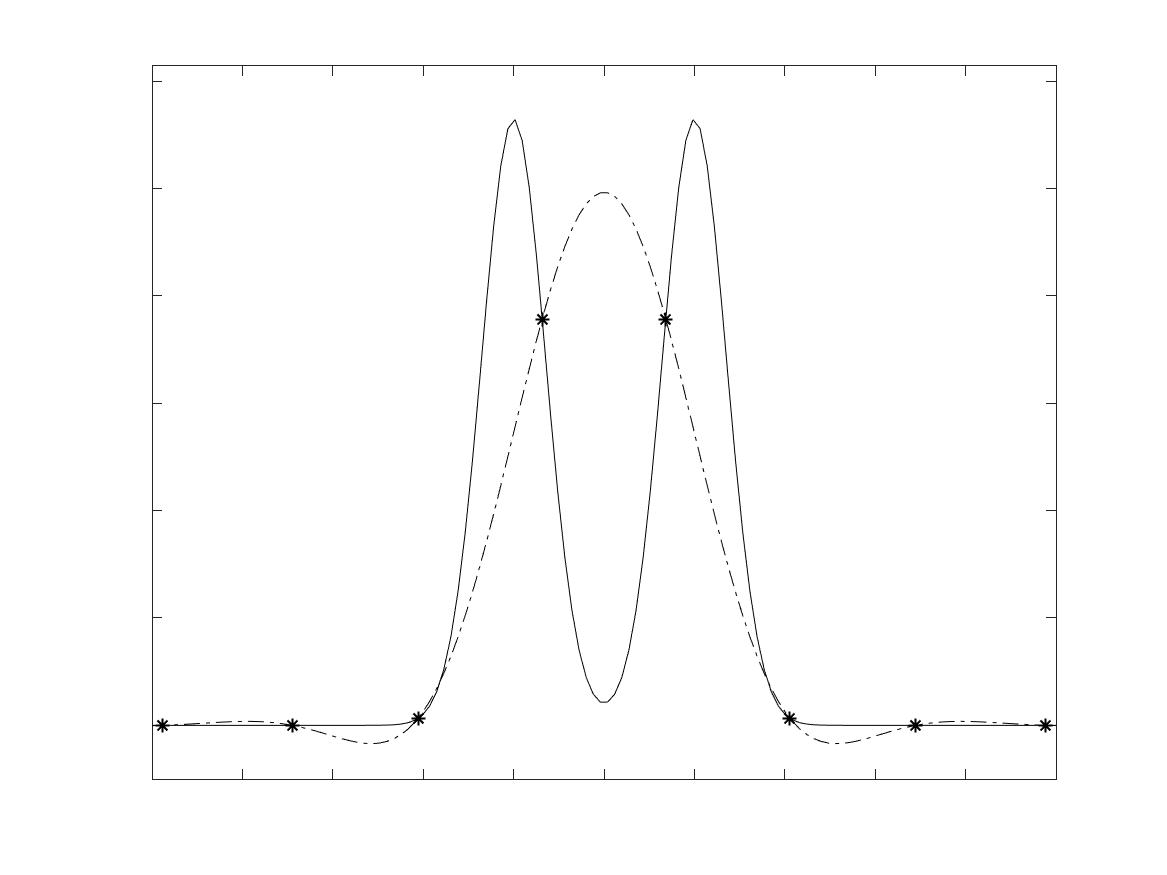}%{./figurePaper/HermiteN2^4alpha04.jpg} 
    \hskip-0.6truecm
     \includegraphics[height=5.5cm]{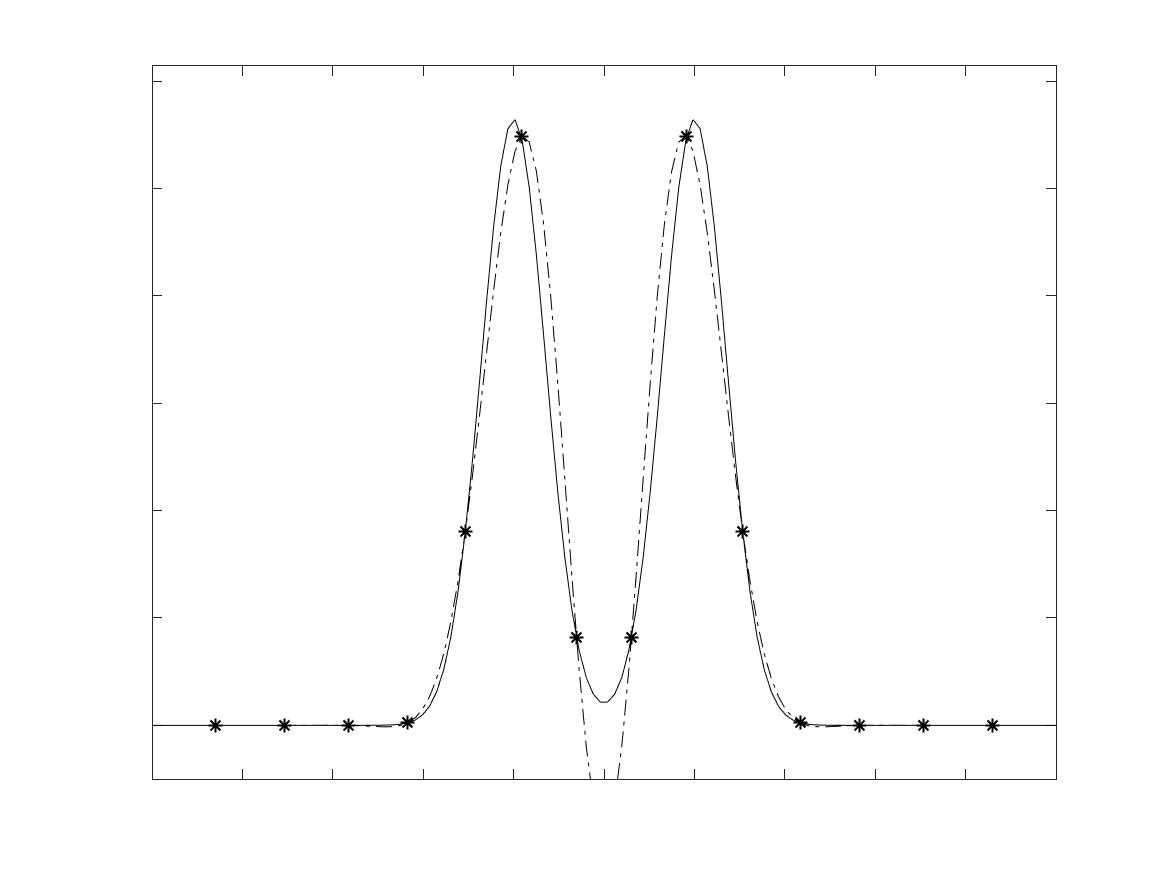}%{./figurePaper/HermiteN2^4alpha09.jpg} 
  }
  \centerline{
  \includegraphics[height=5.5cm]{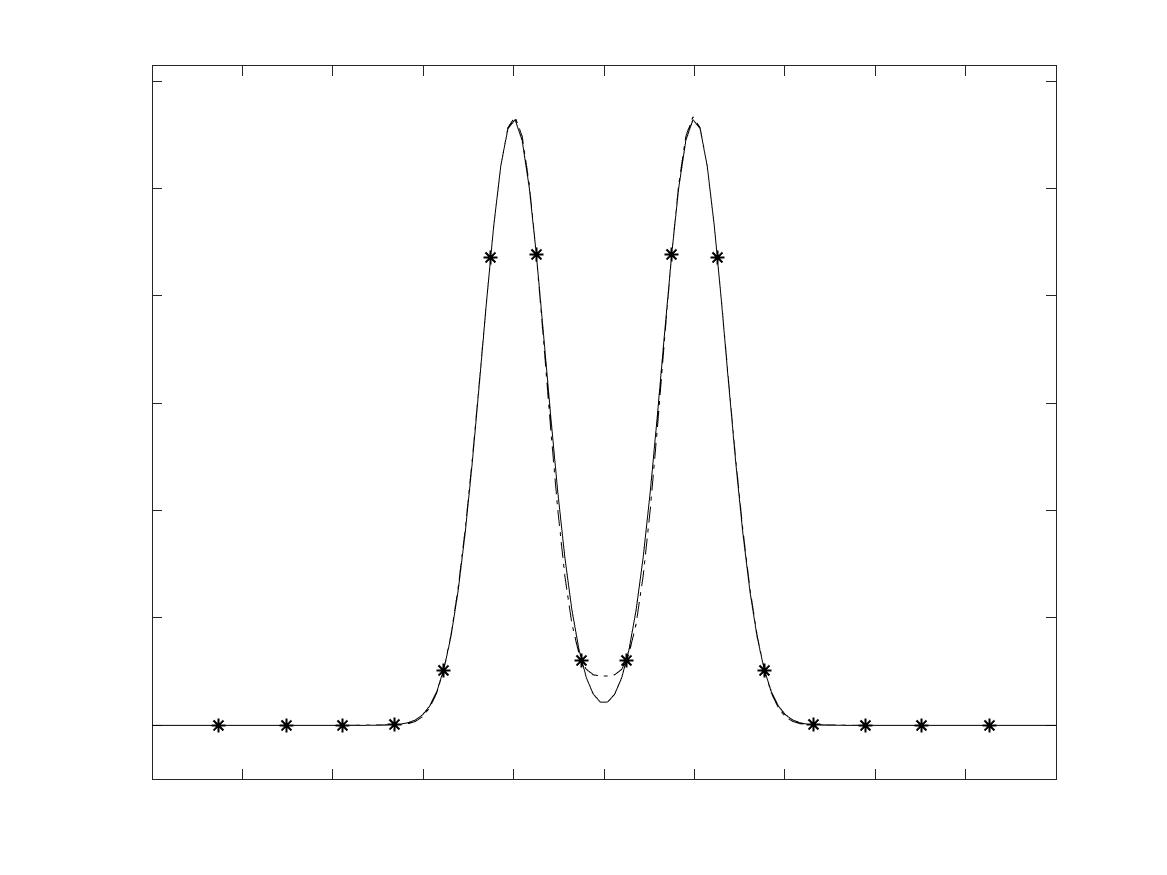}%{./figurePaper/HermiteN2^4alpha11.jpg} 
   \hskip-0.6truecm
    \includegraphics[height=5.5cm]{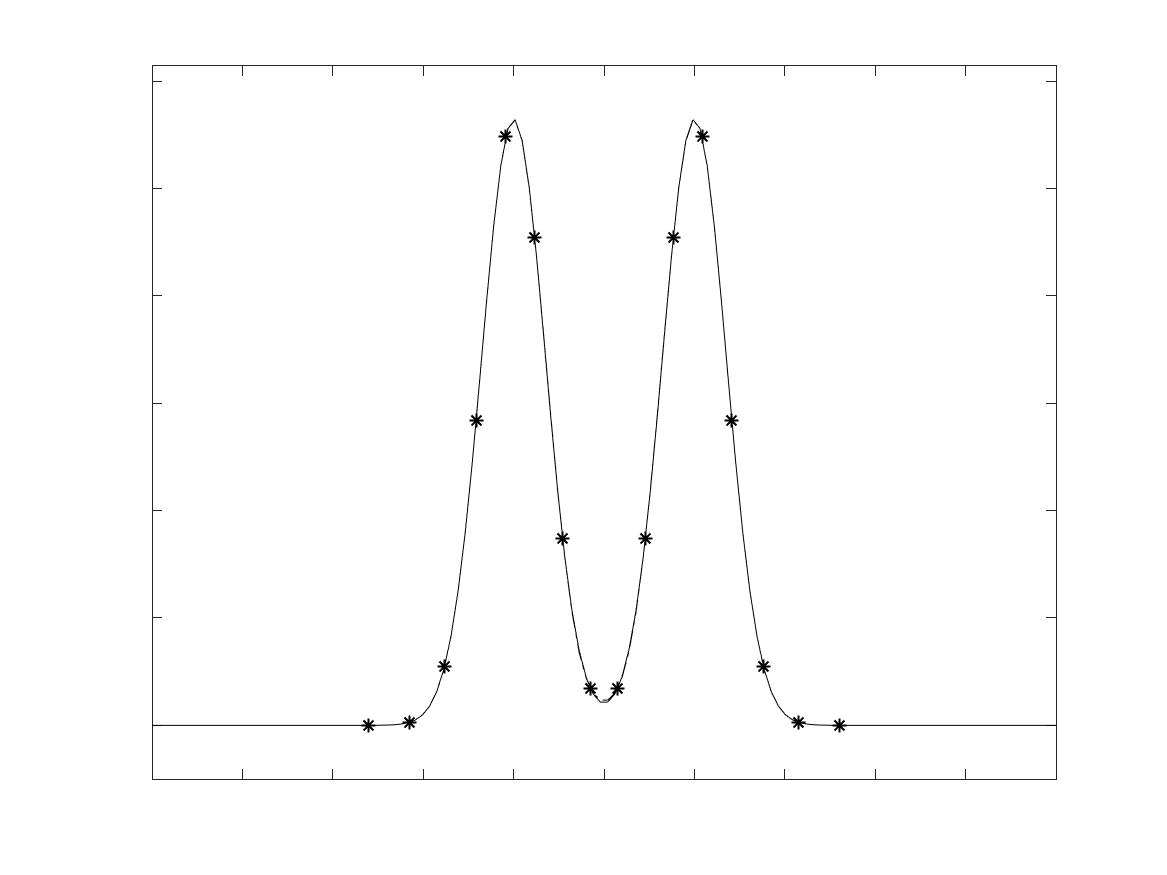}%{./figurePaper/HermiteN2^4alpha18.jpg} 
  }
  \caption{Interpolation with respect to $v\in [-5,5]$ of the function
    \eqref{eq:1D1V:twostreamfzero} with $M=2^{4}$.
    The interpolation is made using  scaled 
     Hermite  nodes with 
    $\alpha=0.4$  (top-left), $\alpha=0.9$  (top-right),
    $\alpha=1.1$ (bottom-left) and $\alpha=1.8$ (bottom-right).  }
  \label{fig2approx}
\end{figure}

Note also that in the Hermite case, the interpolation procedure is very sensitive to the choice
of the parameter $\alpha$.
In Figure \ref{fig2approx} we plot the interpolants for $M=2^{4}$ and
different values of $\alpha$, namely::   $\alpha=0.4$, $\alpha=0.9$,
    $\alpha=1.1$  and $\alpha=1.8$.

%%%%%%%%%%%%%%%%%%%%%%%%%%%%%%%%%%%%%%%%%%%%%%%%%%%%%%%%%%
\section{Conservation Properties}
\label{sec:conservation:properties}
%%%%%%%%%%%%%%%%%%%%%%%%%%%%%%%%%%%%%%%%%%%%%%%%%%%%%%%%%%

%\textcolor{red}{}

The discrete counterpart of \eqref{eq:1D1V:massConserv} (which is the
law stating that the number of particles is preserved) holds for the
scheme \eqref{eq:1D1V:favanza}-\eqref{eq:1D1V:favanzadef}, when $g=0$.
Of course, the conservation of this basic quantity is rather important
from the physics viewpoint.
At  step $t^{k}=k \Delta t$, we define:
\begin{equation}\label{eq:1D1V:massadiscreta}
  \Qs^{(k)}_{N,M} 
  = \frac{2\pi}{N}\sum_{n=0}^{N-1}\sum_{m=0}^{M-1}c^{(k)}_{nm} w_m
  = \int_\Omega f^{(k)}_{N,M}(x,v) dxdv \approx \int_\Omega f(t^k,x,v)dxdv,
\end{equation}
where $w_m=2\pi/M$ in the periodic Fourier case (recall the quadrature
formula \eqref{eq:1D1V:fdq}), while the weights are related to the
quadrature formula \eqref{eq:1D1V:fdqleg}, which are also valid for both Legendre and 
Hermite cases.
%%
%The correspondence of the two integrals in
%\eqref{eq:1D1V:massadiscreta} is true up to an error that is
%spectrally accurate, due to the excellent properties of Gaussian
%quadrature.  
%% 
From \eqref{eq:1D1V:favanza}, at time step $t^{k+1}$ we find out
that:
\begin{align}
  \Qs^{(k+1)}_{N,M}
  = \frac{2\pi}{N}\sum_{n=0}^{N-1}\sum_{m=0}^{M-1}\Big(\cs_{nm}^{(k)} + \Delta t\,\Phi_{nm}^{(k)}\Big) w_m
  = \Qs^{(k)}_{N,M}+\Delta \Qs^{(k)}_{N,M}= \Qs^{(k)}_{N,M}.
\end{align}
In fact, one has:
\begin{align}
  \Delta \Qs^{(k)}_{N,M} =
  -\Delta t \sum_{m=0}^{M-1} v_m\left[ \int_{\Omega_x}\frac{\partial f^{(k)}_{N,M}}{\partial x}
    (x, v_m)dx\right] \frac{w_m}{\omega (v_m)} \nonumber\\
  + \ \Delta t \frac{2\pi}{N}\sum_{n=0}^{N-1} E^{(k)}_N(x_n)\left[ \int_{\Omega_v}\frac{\partial f^{(k)}_{N,M}}{\partial v}(x_n, v)dv\right] = 0,
  \label{eq:1D1V:massadiscretacons}
\end{align}
where $\omega (v)=1$ in the trigonometric and in the Legendre case,
while $\omega (v)=\exp (-v^2)$ in the Hermite case.
The integrals above are zero as a consequence of the boundary
conditions.
This is true for all the cases we are considering in this paper, i.e.,
periodic boundary conditions or homogeneous Dirichlet conditions
(imposed at the points $\pm 1$ when $\Omega_v=]-1,1[$ or through a
suitable exponential decay when $\Omega_v=\R$).
In the end, we proved that the quantity in
\eqref{eq:1D1V:massadiscreta} does not change passing from $k$ to
$k+1$.
The same property holds for the scheme \eqref{eq:1D1V:favanzabdf2co}.
The proof follows after recognizing that, for $g=0$, one has:
\begin{align}
  \Qs^{(k+1)}_{N,M} = \frac{4}{3} \Qs^{(k)}_{N,M} - \frac{1}{3} \Qs^{(k-1)}_{N,M}.
\end{align}
The above recurrence relation provides constant values when the initial coefficients
are such that $c^{(1)}_{ij} = c^{(0)}_{ij}$.

\smallskip
Analogous considerations are made regarding the conservation in time
of the momenta $\int_\Omega v^r f(t,x,v) dxdv$, where $r\geq 0$ is an
integer.
In the discrete case we define at time $t^{k}$, $k=0,1,\ldots,K$:
\begin{equation}\label{eq:1D1V:momentodiscreto}
  Q^{(k)}_{N,M,r} 
  = \frac{2\pi}{N}\sum_{n=0}^{N-1}\sum_{m=0}^{M-1}v_{m}^r c^{(k)}_{nm} w_m
  \approx \int_\Omega v^r f^{(k)}_{N,M}(x,v) dxdv.
\end{equation}
Of course $Q^{(k)}_{N,M,0}=Q^{(k)}_{N,M}$.
We expect conservation of the discrete momenta in the Legendre or
Hermite case, for values of $r$ up to the integration capabilities of
Gaussian formulas.
This means $r\leq M$ for the Legendre case and $r\leq M+2$ for the
Hermite case.
We note instead that, being $v$ not a trigonometric function, the
integration formula does not hold in the periodic case.
In this situation, $v$ can be substituted by its projection in the
$L^2(\Omega)$ norm on the finite dimensional space ${\bf Y}_{N,M}$, up
to an error that decays spectrally.
This procedure may however generate a Gibb's phenomenon across the
points of $\Omega$ with $v=2\pi$, where $v^rf$ is discontinuous.
Nevertheless, if the function $f$ shows a fast decay near the boundary with
respect to the variable $v$, the conservation of momenta can be
achieved up to negligible errors.

\smallskip
At time $t^{k}$, $k=0,1,\ldots,K$, we can also consider the discrete
version of \eqref{eq:1D1V:totEnergy} i.e.:
\begin{equation}\label{eq:1D1V:energiadiscreta}
  {\cal E}(t^k) \approx {\cal E}^{(k)}_{N,M}
  = \frac12\left( 
    \frac{2\pi}{N}\sum_{n=0}^{N-1}\sum_{m=0}^{M-1}v^2_mc^{(k)}_{nm} w_m+
    \frac{2\pi}{N}\sum_{n=0}^{N-1}\Big[ E^{(k)}_N (x_n)\Big]^2
  \right).
\end{equation}
that includes the term $Q^{(k)}_{N,M,2}$.
Conservation properties of this quantity will be tested later in the
numerical section.

\raggedbottom
%%%%%%%%%%%%%%%%%%%%%%%%%%%%%%%%%%%%%%%%%%%%%%%%%%%%%%%%%%
\section{Numerical experiments}
\label{sec:numerical:experiments}
%%%%%%%%%%%%%%%%%%%%%%%%%%%%%%%%%%%%%%%%%%%%%%%%%%%%%%%%%%

Together with other notable examples, the numerical scheme here proposed has been already validated 
in ~\cite{Fatone-Funaro-Manzini:2018} in the case of the two-stream instability example,
which is a standard benchmark for plasma physics codes.
This corresponds to the setting $\Omega_{x}=[0,4\pi[$,
$\Omega_{v}=[-5,5]$ in \eqref{eq:1D1V:V}, \eqref{eq:1D1V:Vci},
\eqref{eq:1D1V:EMaxwell}, \eqref{eq:1D1V:charge-density},
with the initial guess  given by:
\begin{align} \label{eq:1D1V:twostreamf}
  \bar{\fs}(\xs,\vs)=\frac{1}{2a\sqrt{2\pi}}  
  \left[ 
    \exp { \left(-\frac{\vs-\beta}{a\sqrt{2}}\right)^{2}}+
    \exp{ \left(-\frac{\vs+\beta}{a\sqrt{2}}\right)^{2}}
  \right]
  \left[ 
    1 + \epsilon \cos\left( \kappa \xs \right)
  \right],
  %% \qquad\xs\in\Omega_{x},\,\,\vs\in\Omega_{v},
\end{align}
with $a=1\slash{\sqrt 8}$, $\beta=1$, $\epsilon=10^{-3}$,
$\kappa=0.5$.
The exact solution is approximated by periodic basis functions in the
variable $\xs$ and Fourier, Legendre  and Hermite basis functions in the variable
$\vs$. For practical purposes, depending the basis examined, we rewrite the differential equation on the domain
$\Omega =  [0,2\pi [\times \Omega_v$, with $\Omega_v =[0, 2\pi[$ (Fourier) or
$\Omega_v =]-1,1[$ (Legendre), while $\Omega_v =\R$ is not modified in the Hermite case. The purpose of the paper is to discuss the different performances in dependence of 
the basis used for the variable $v$. Comparisons between the various approaches 
are made by taking less degrees of freedom than those actually necessary to resolve accurately the equation.
In particular, we will perform experiments by setting $N=M=2^{4}$. 

\smallskip

\begin{figure}
  \centerline{
   \includegraphics[height=5.5cm]{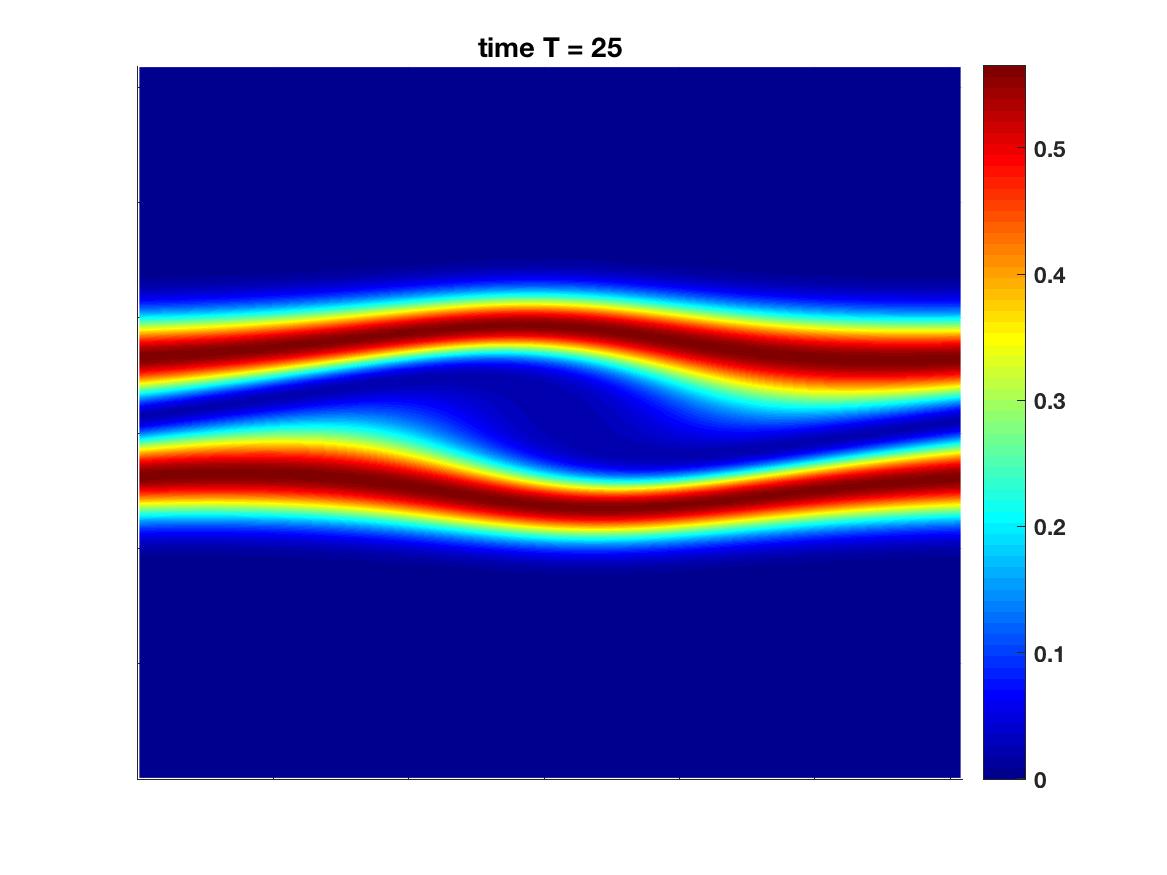}%{./figurePaper/FirstvorticeN25M26A.jpg} 
     \hskip-0.6truecm  
      \includegraphics[height=5.5cm]{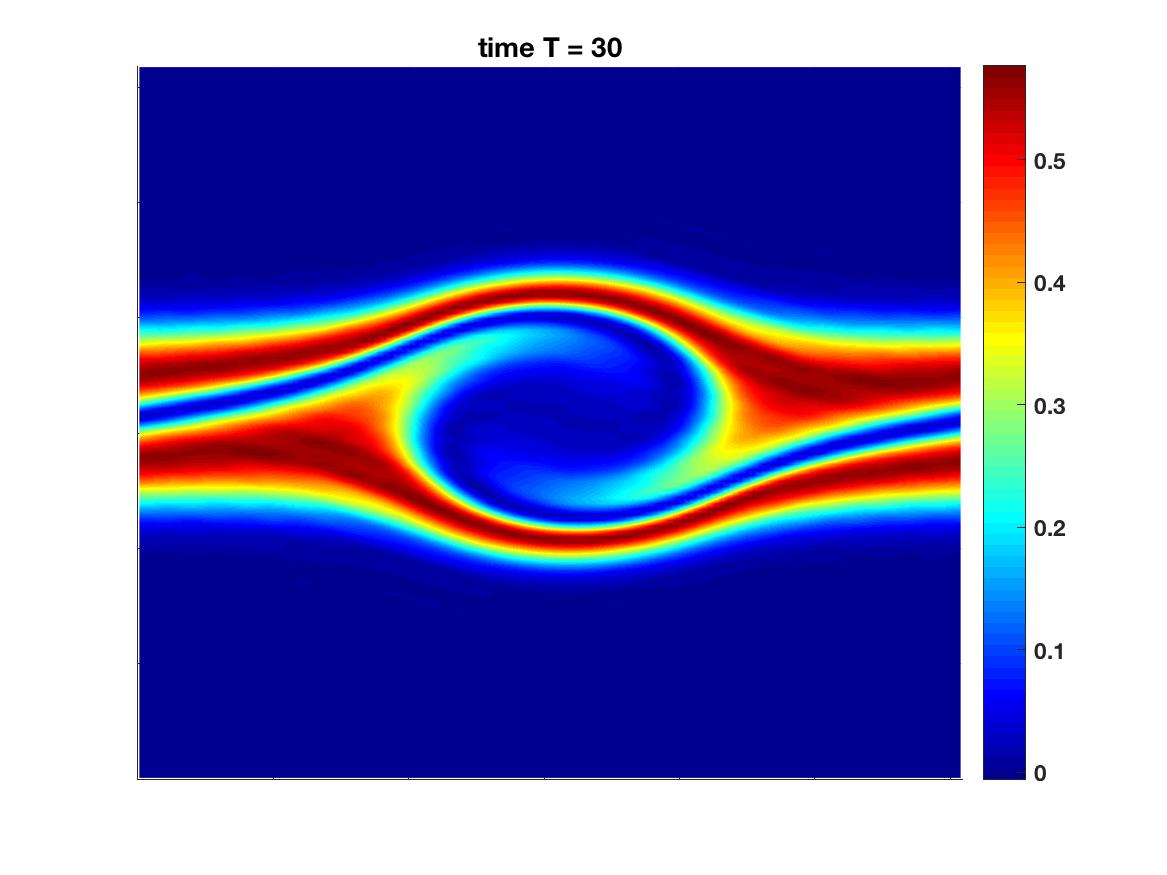}%{./figurePaper/FirstvorticeN25M26B.jpg} 
  }
  \caption{ Approximated distribution function obtained  by using the Fourier-Fourier method with  $N=2^{5}$, $M=2^{6}$ and  $\Delta t=1/800$.   }
   \label{fig3}
%\end{figure}
%
\vskip2truecm
%
%\begin{figure}
  \centerline{
   \includegraphics[height=5.5cm]{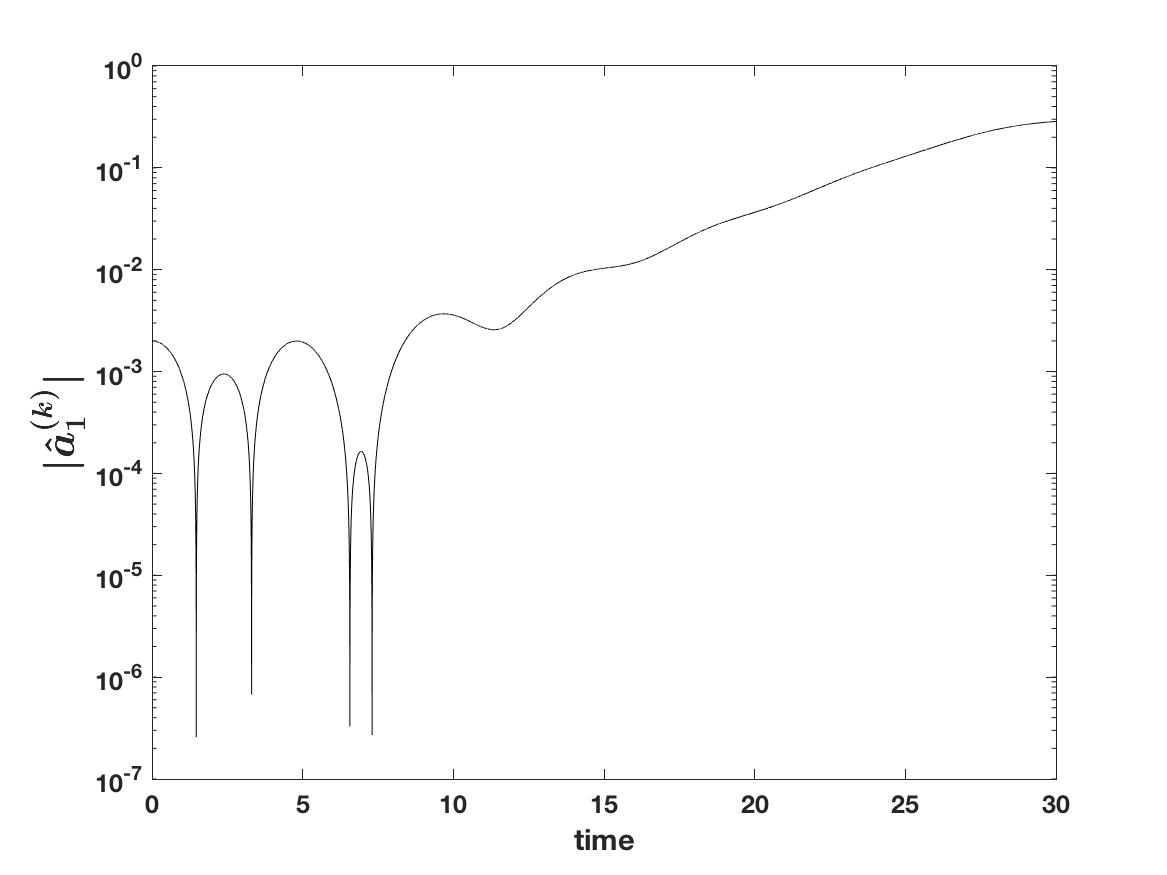}%{./figurePaper/FirstvorticeN25M26C.jpg} 
  }
  \caption{ Time evolution of the first
Fourier mode of the electric field $\Es^{(k)}_N$, i.e.
$|\hat{\as}_{1}^{(k)}|$ in \eqref{eq:1D1V:FourierSeriesE}  when using  the Fourier-Fourier method for   $N=2^{5}$, $M=2^{6}$ and $\Delta t=1/800$.  }
   \label{fig4}
\end{figure}

In all the experiments to follow, we integrate up to time $T=30$ using the second-order BDF scheme with 
an appropriate time step, in order to guarantee stability. The choice of an accurate time
discretization procedure will allow us to concentrate our attention to the spectral approximation in the variable $x$ and $v$.
First of all, in Figure \ref{fig3}  we show the results at time $T=25$ and time $T=30$
of the solution recovered by the Fourier-Fourier method, by choosing  $N=2^{5}$, $M=2^{6}$ and a rather 
small time step ($\Delta t =0.25 \cdot 10^{-2}$).
This will be the referring figure for the successive comparisons.  

% since the CFL
% condition \eqref{eq:1D1V:CFL} requires $\Delta t$ to be proportional
% to $1/\max\{N,M\}$.
%%

% A relevant quantity to be examined is $|\hat{\as}_{1}^{(k)}|$, which is the first
% Fourier mode of the electric field $\Es^{(k)}_N$ in \eqref{eq:1D1V:FourierSeriesE}.
% %%
% Figure \ref{fig4} shows the time evolution of such a coefficient when
% Fourier-Fourier scheme with  $N=2^{5}$, $M=2^{6}$ and  $\Delta t=0.0025/2$ is used.
% The behavior is exactly the one expected.
%%

\smallskip

Figure \ref{fig4} shows the time evolution of $|\hat{\as}_{1}^{(k)}|$,
which is the first Fourier mode of the electric field $\Es^{(k)}_N$ in
\eqref{eq:1D1V:FourierSeriesE}, when the Fourier-Fourier method with
$N=2^{5}$, $M=2^{6}$ and $\Delta t=1/800$ is used. 
This crucial quantity is often displayed in order to make comparisons
between computational codes, also because its behavior is predicted by
theoretical consideration.
Indeed, the slope of the segment starting from $T=15$ in the plot of Figure
\ref{fig4} agrees with the expectancy~\cite[Chapter~5]{Bittencourt:2004}.

% {\textcolor{red}{Qualcosina in piu' sul comportamento teorico,
%     references}}.  The plot given in actually reflects the expected
% situation.
%%
% Note that the use of the second-order BDF scheme with the same
% parameters furnishes results indistinguishable to the eye respect to
% those shown in Figures \ref{fig3}, \ref{fig4}.

% {\textcolor{red}{Non so se convenga rifare Figures \ref{fig3},
%     \ref{fig4} con il BDF2 visto che poi mostriamo (quasi)
%     sempre %risultati ottenuti con tale metodo.
%     Forse conviene anche motivare la scelta del BDF2 (propriet\`a
%     spettrali)}}

% This will be the referring figure for the first Fourier mode of the
% electric field $\Es^{(k)}_N$.
\smallskip

\begin{figure}
  \centerline{
    \includegraphics[height=5.5cm]{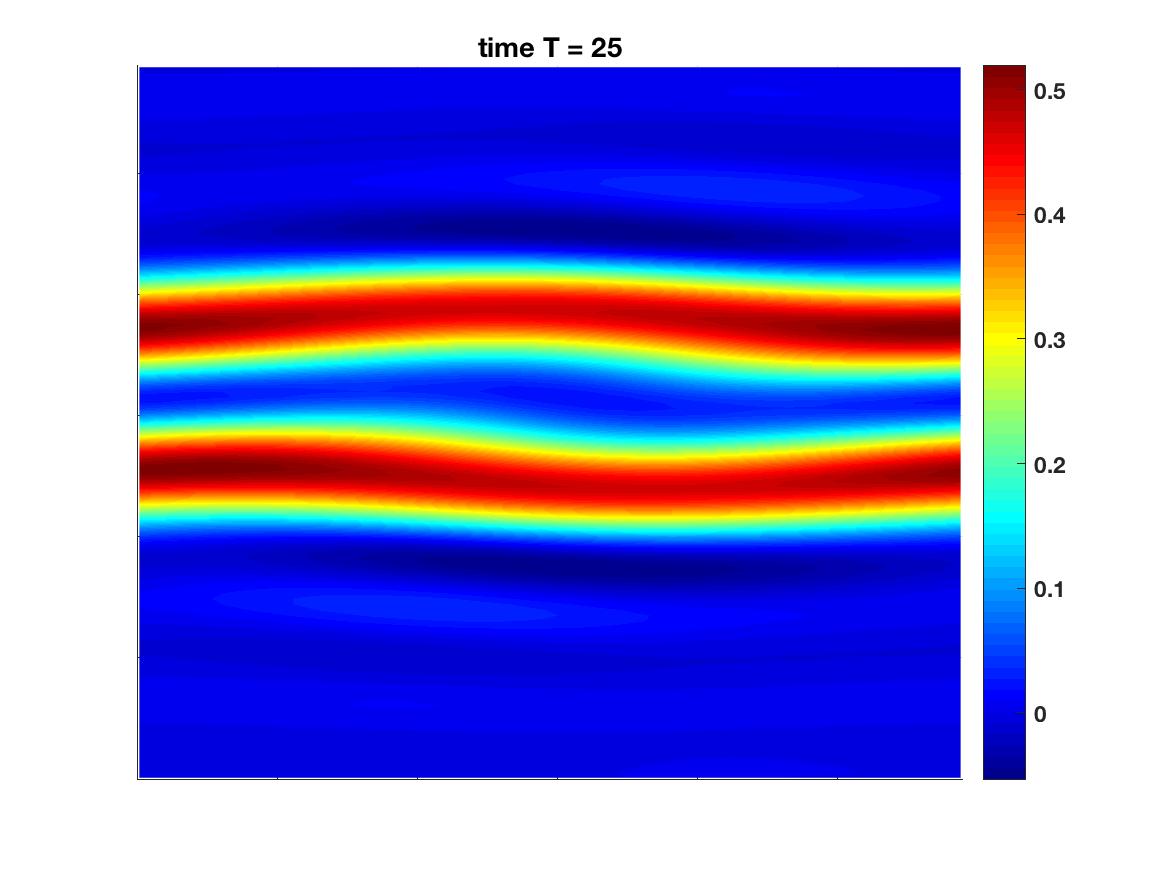}%{./figurePaper/aFFN16BDF2.jpg} 
     \hskip-0.6truecm
    \includegraphics[height=5.5cm]{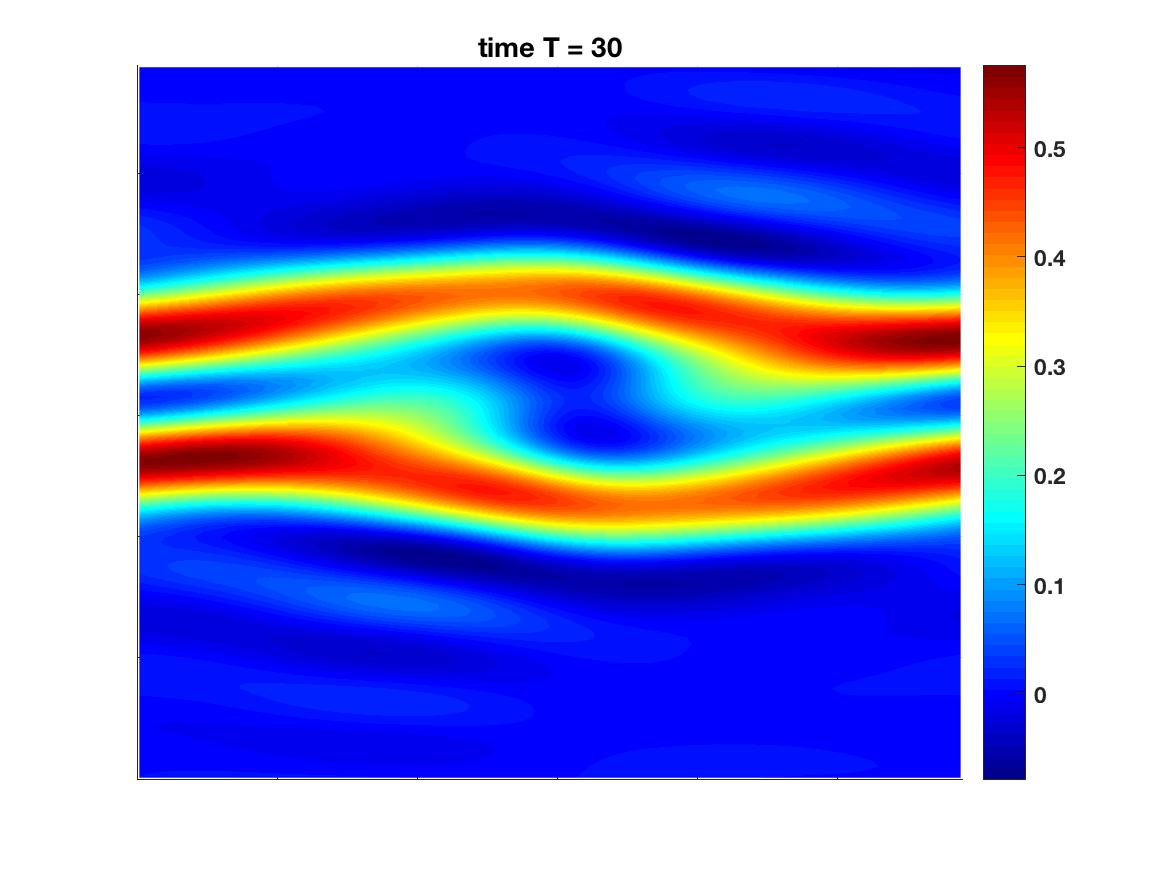}%{./figurePaper/bFFN16BDF2.jpg} 
  }
  \caption{ Approximated distribution function obtained  by using  the Fourier-Fourier method   with    $N=M=2^{4}$ and  $\Delta t=0.01$.}
   \label{fig5}
\end{figure}

\begin{figure}
  \centerline{
    \includegraphics[height=5.5cm]{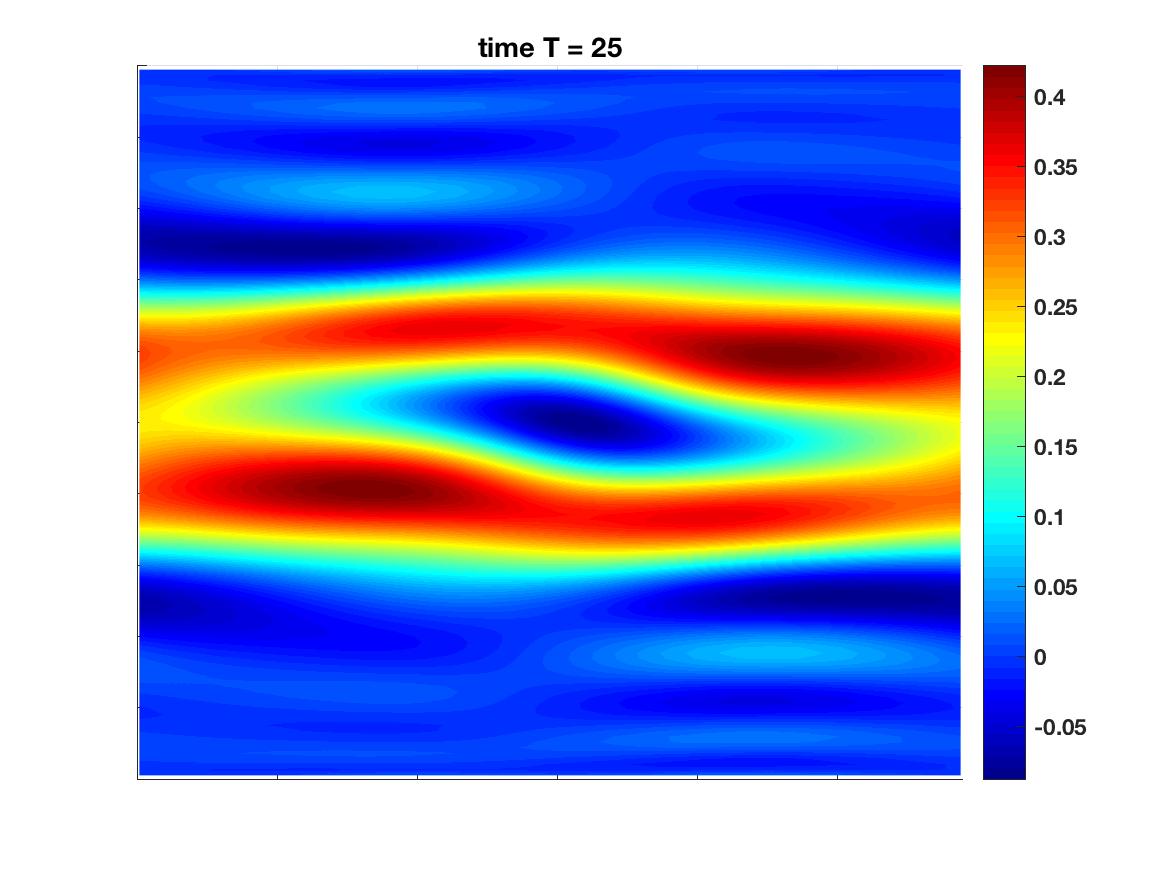}%{./figurePaper/aFLN16BDF2.jpg} 
  }
  \caption{ Approximated distribution function obtained  by using  the Fourier-Legendre method   with    $N=M=2^{4}$ and $\Delta t=0.01$.}
   \label{fig6}
\end{figure}

The plots of Figure \ref{fig5} show the numerical distribution at time
$T=25$ and time $T=30$ obtained by using the Fourier-Fourier method
with $N=M=2^{4}$ and $\Delta t=0.01$.
Such an approximation is rough, but the computation is stable and
provides an idea of the real behavior.
Unfortunately, as we predicted in Section~\ref{sec:basis:choice},
results of this kind are not available for the Fourier-Legendre case.
The algorithm is stable up to $T=25$ (see Figure \ref{fig6}), but then it explodes. 
As we already mentioned, the explanation of this fact may rely on the
poor approximation of the Legendre interpolants at the interior of the
interval.
The performances improve relevantly by increasing the value of
$M$. The Fourier-Fourier method looks however a better choice, by
virtue also of the fact that it allows for the use of the DFT.

\smallskip

\begin{figure}
   \centerline{
    \includegraphics[height=5.5cm]{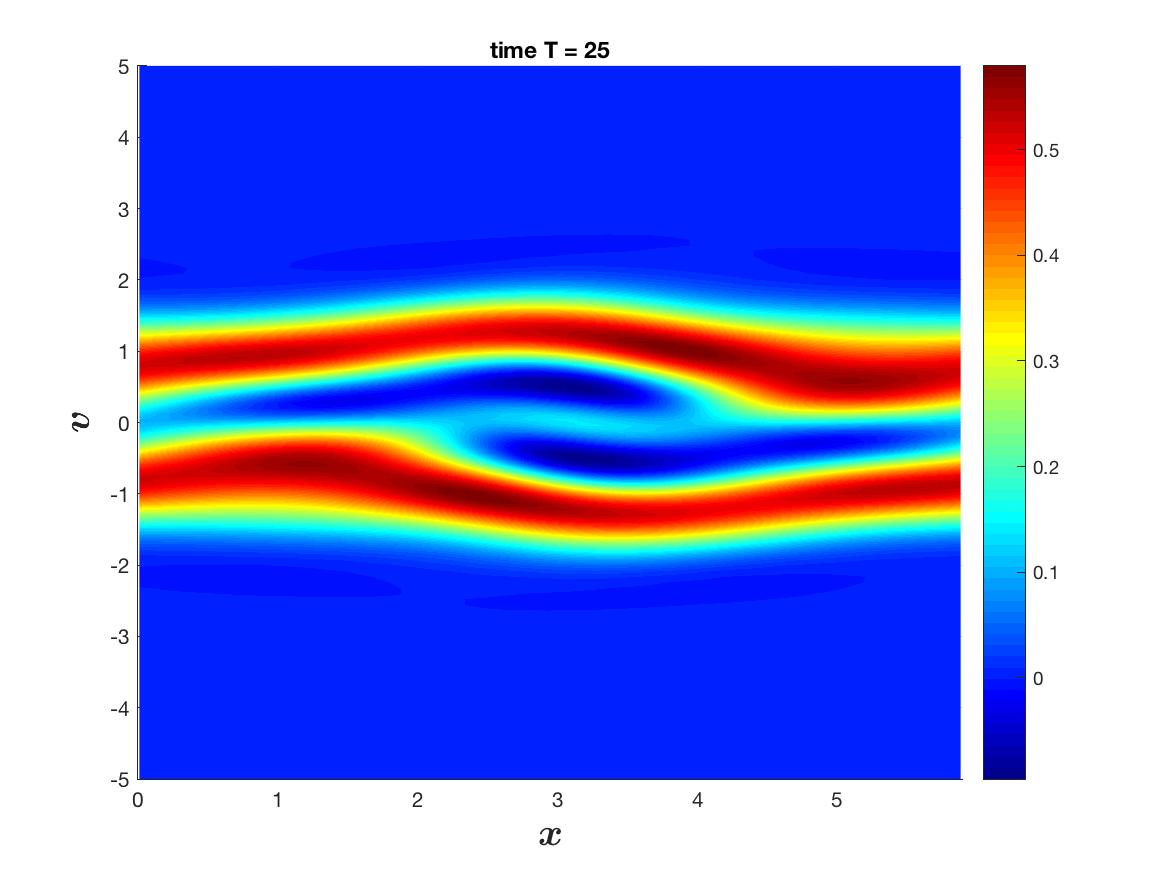}%{./figurePaper/aFHN16alpha1BDF2.jpg} 
    \hskip-0.6truecm
    \includegraphics[height=5.5cm]{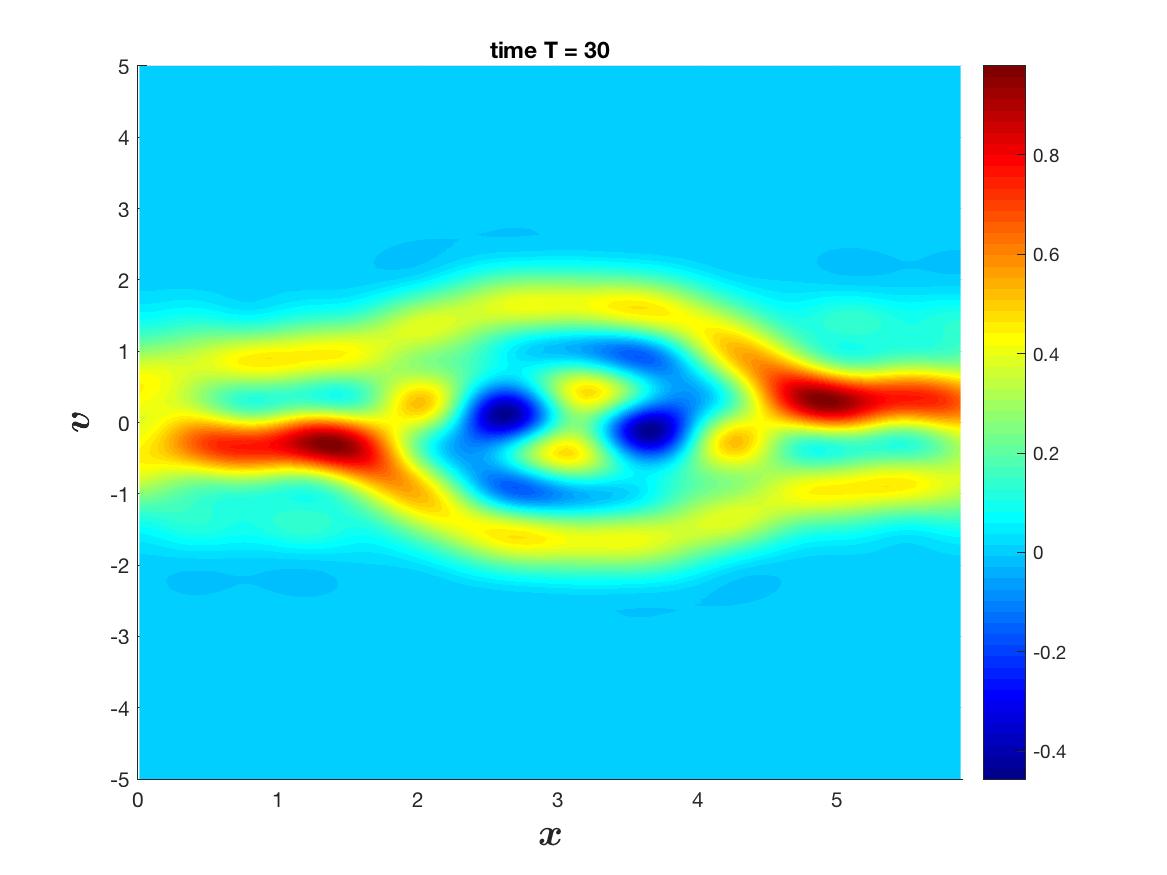}%{./figurePaper/bFHN16alpha1BDF2.jpg} 
  }
  \caption{ Approximated distribution function obtained  by using  the  Fourier-Hermite  method   with    $N=M=2^{4}$, $\Delta t=0.01$ and $\alpha=1$ .}
   \label{fig7}
\end{figure}

\begin{figure}
  \centerline{
    \includegraphics[height=5.5cm]{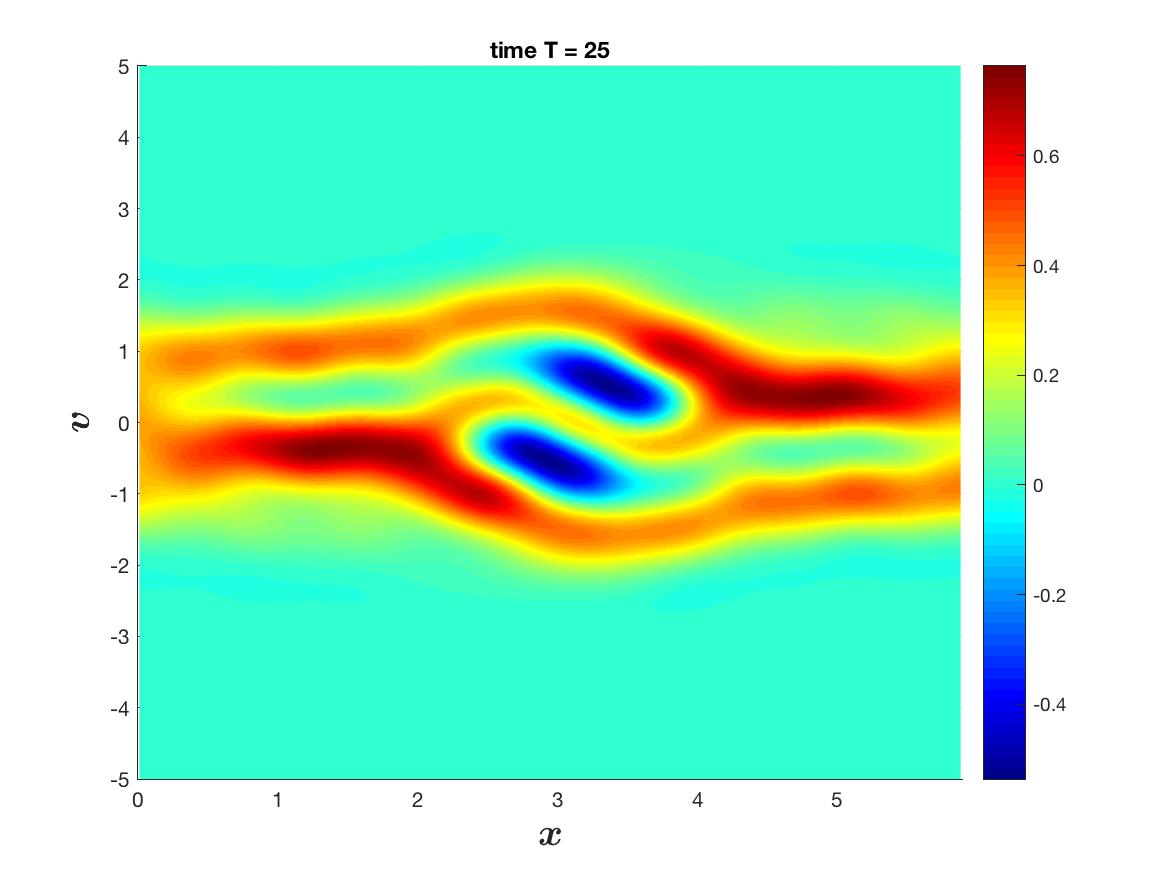}%{./figurePaper/aFHN16alpha09BDF2.jpg} 
    \hskip-0.6truecm
    \includegraphics[height=5.5cm]{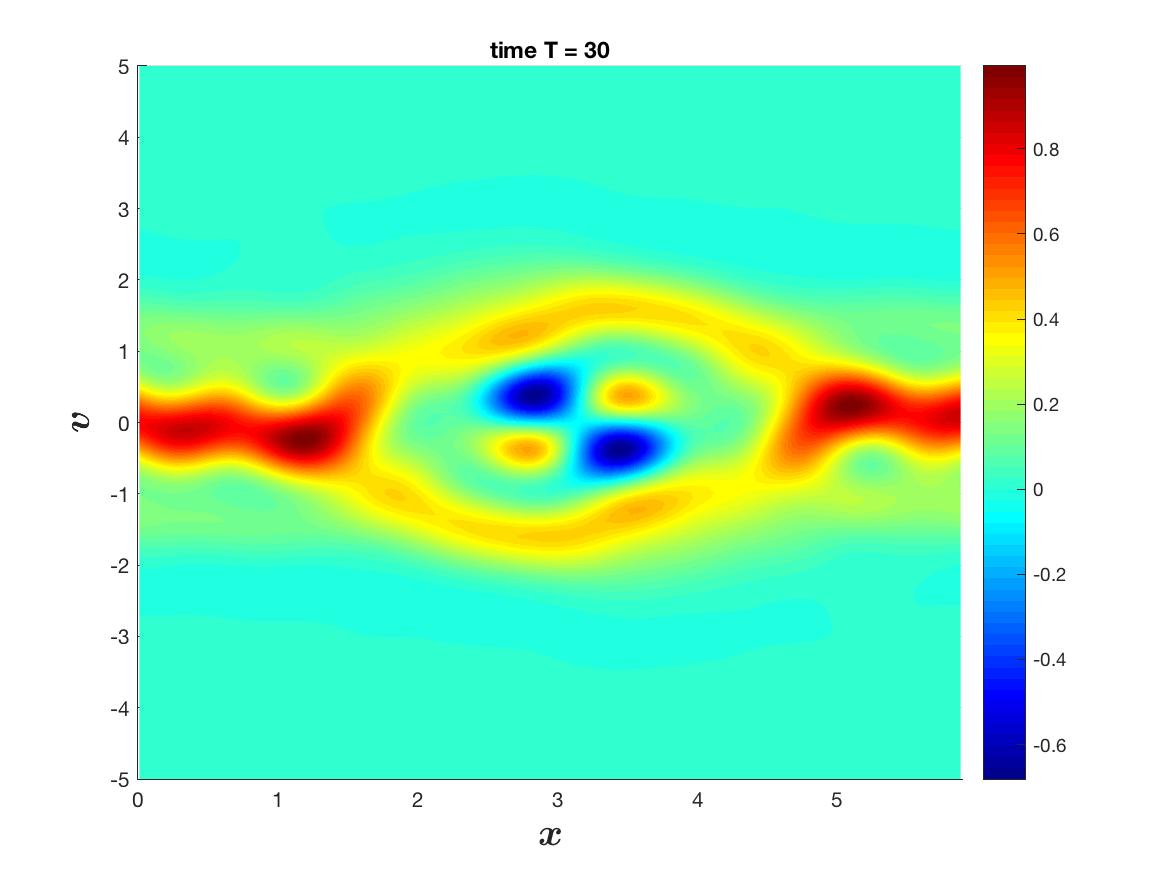}%{./figurePaper/bFHN16alpha09BDF2.jpg} 
  }
  \centerline{
    \includegraphics[height=5.5cm]{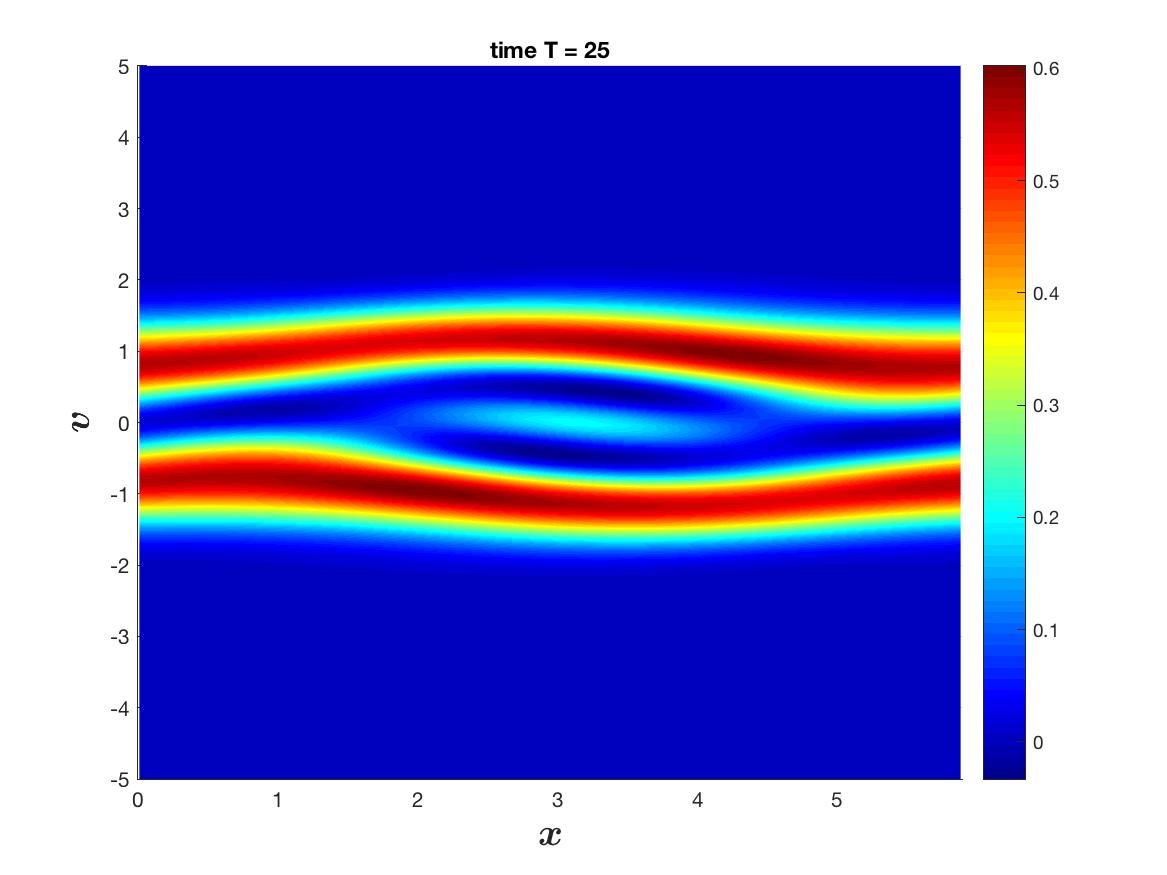}%{./figurePaper/aFHN16alpha11BDF2.jpg}  
     \hskip-0.6truecm   
    \includegraphics[height=5.5cm]{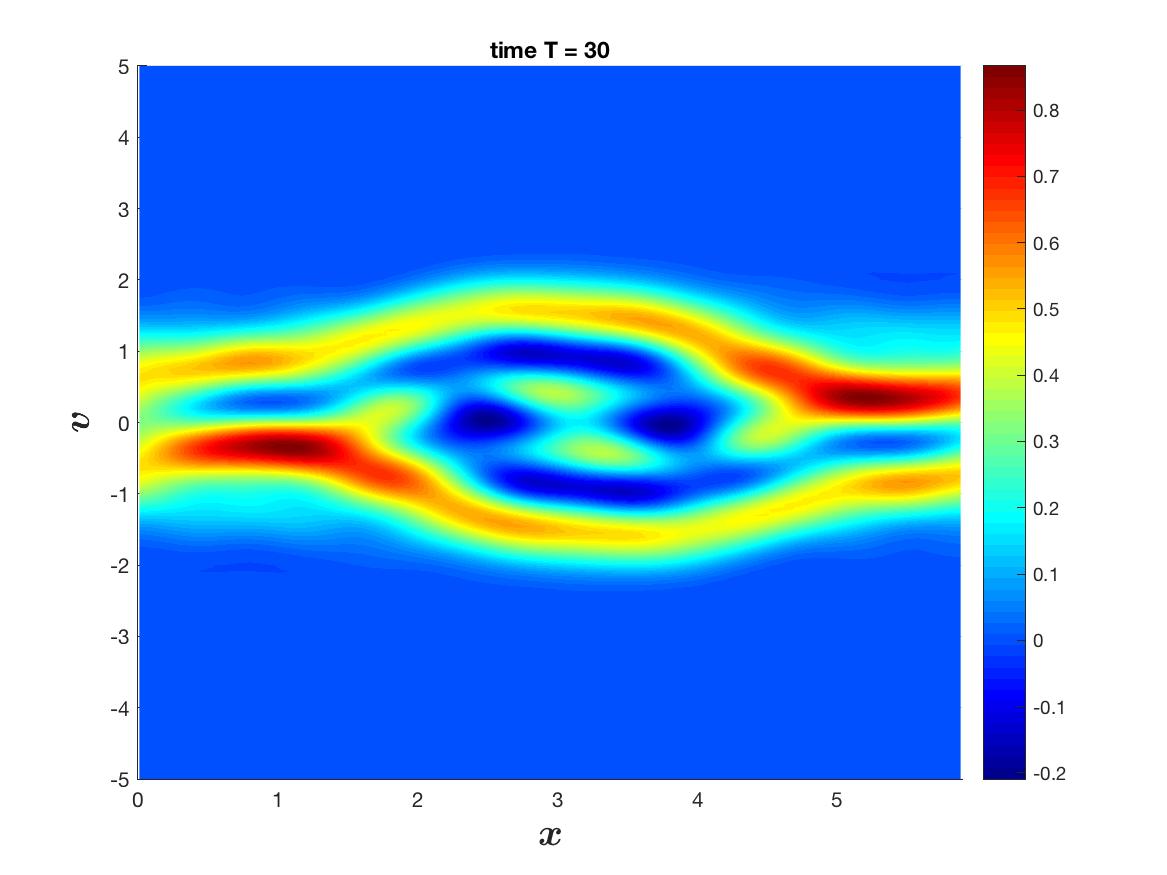}%{./figurePaper/bFHN16alpha11BDF2.jpg}       
  }
  \centerline{
    \includegraphics[height=5.5cm]{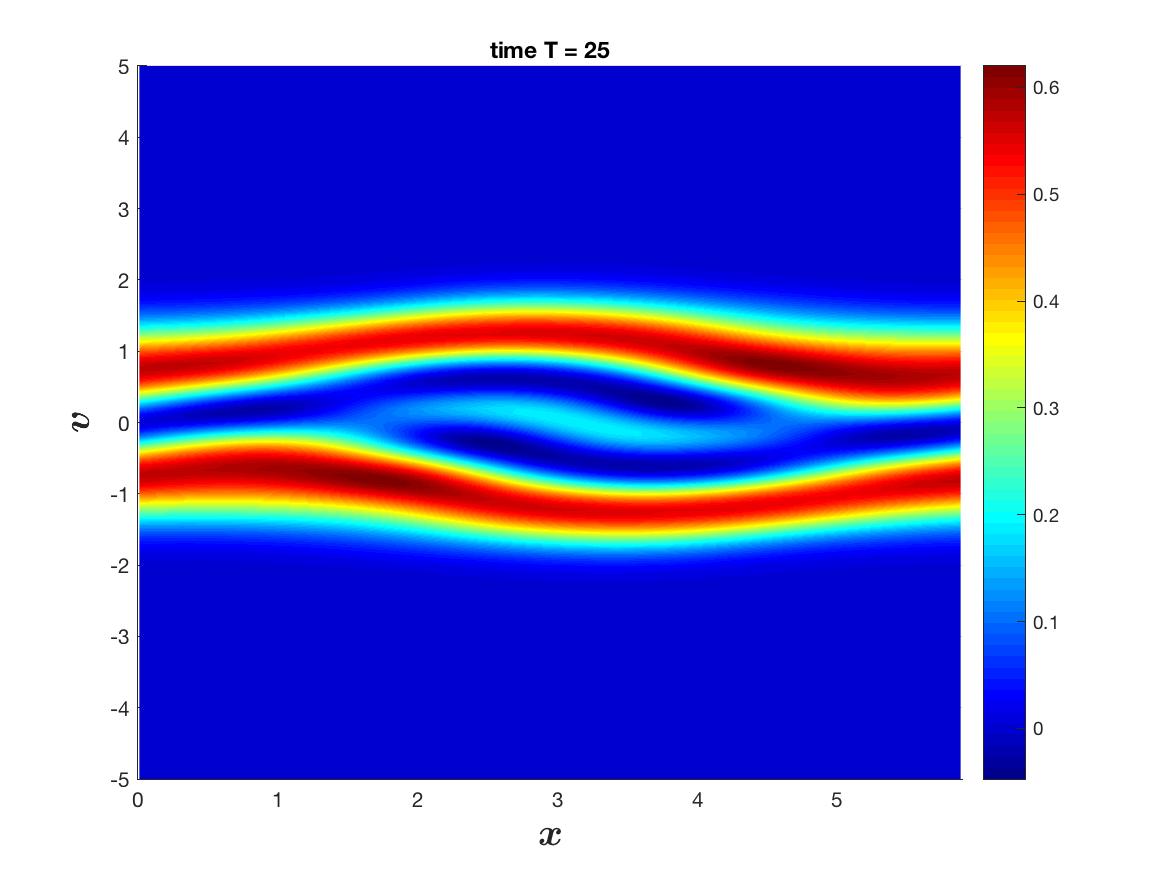}%{./figurePaper/aFHN16alpha12BDF2.jpg}
    \hskip-0.6truecm    
    \includegraphics[height=5.5cm]{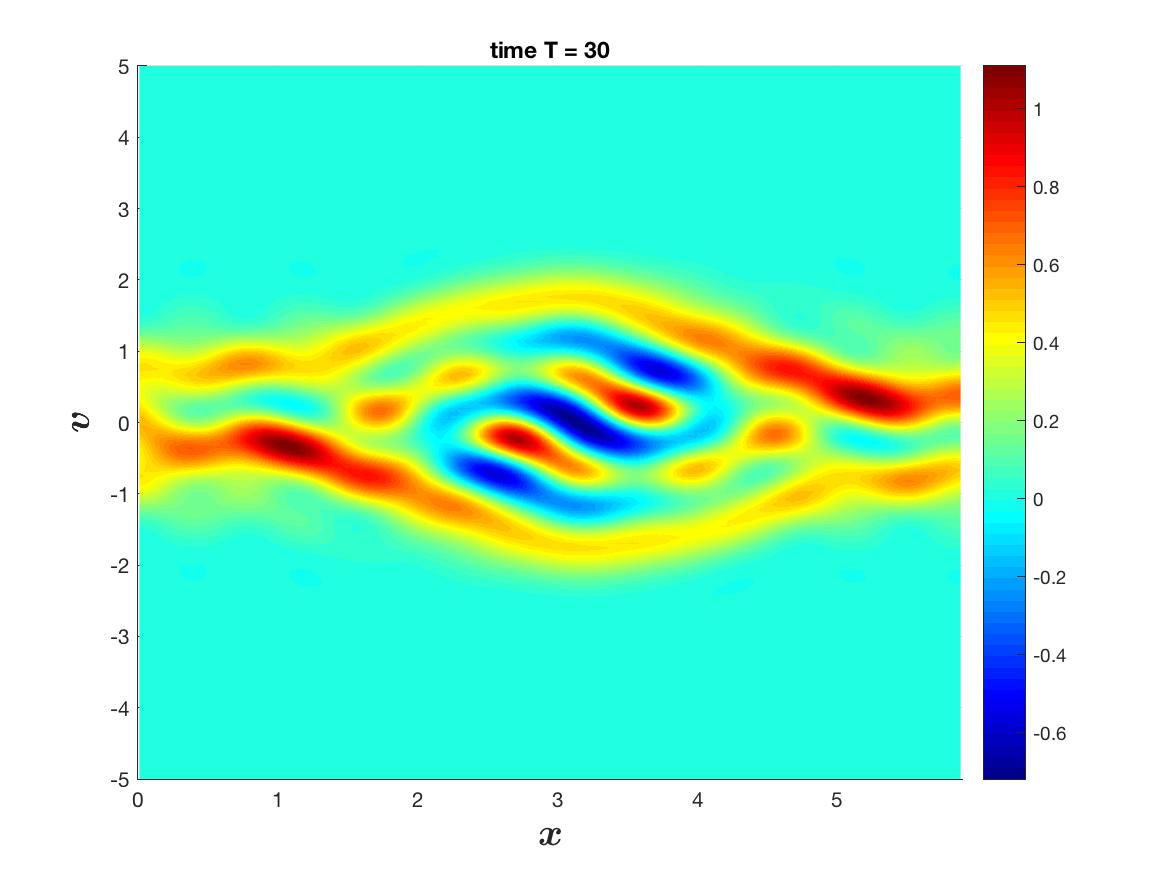}%{./figurePaper/bFHN16alpha12BDF2.jpg}     
  }
  \caption{  Approximated distribution function obtained  by using  the Fourier-Hermite  method  with
    $N=M=2^{4}$, $\Delta t=0.01$ and $\alpha=0.9$ (top), $\alpha=1.1$
    (center), $\alpha=1.2$ (bottom).}
  \label{fig8}
\end{figure}

\begin{figure}
  \centerline{
    \includegraphics[height=5.5cm]{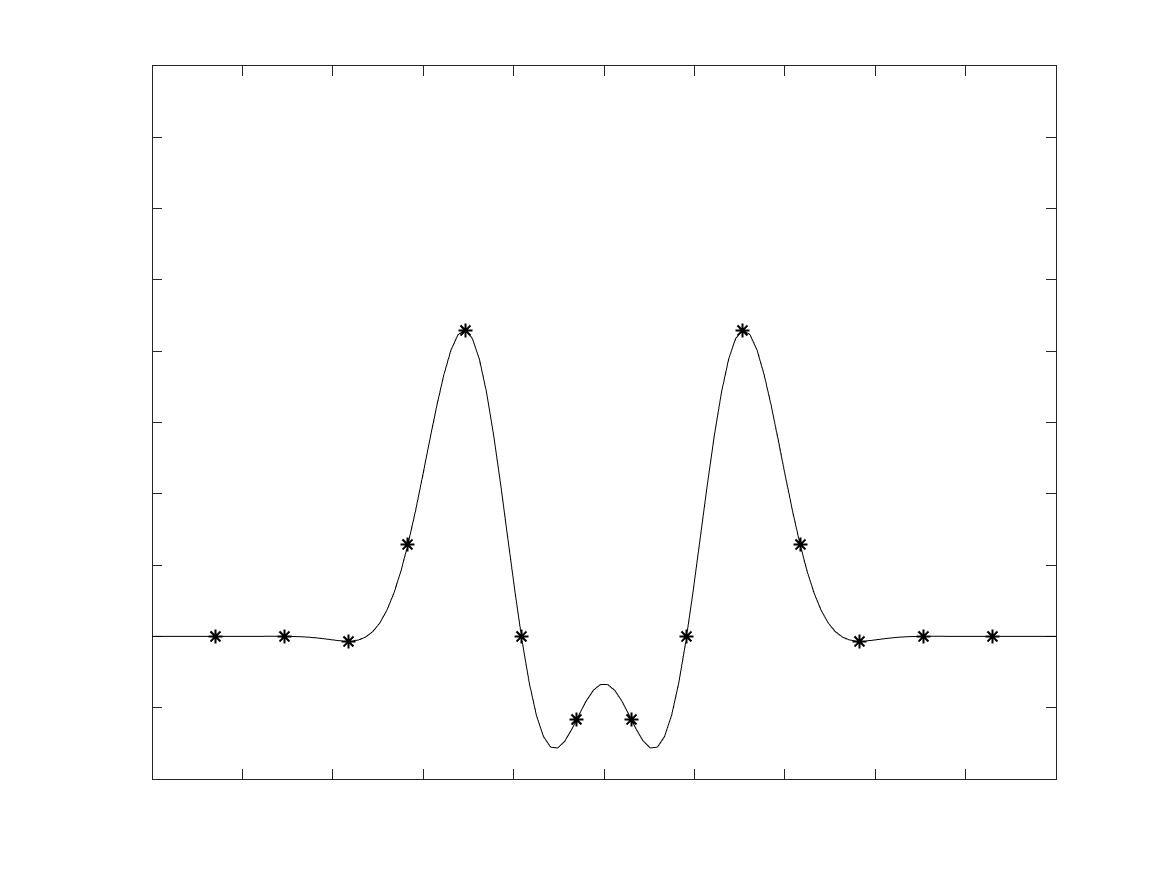}%{./figurePaper/proiezioneT30alpha09.jpg} 
     \hskip-0.6truecm
     \includegraphics[height=5.5cm]{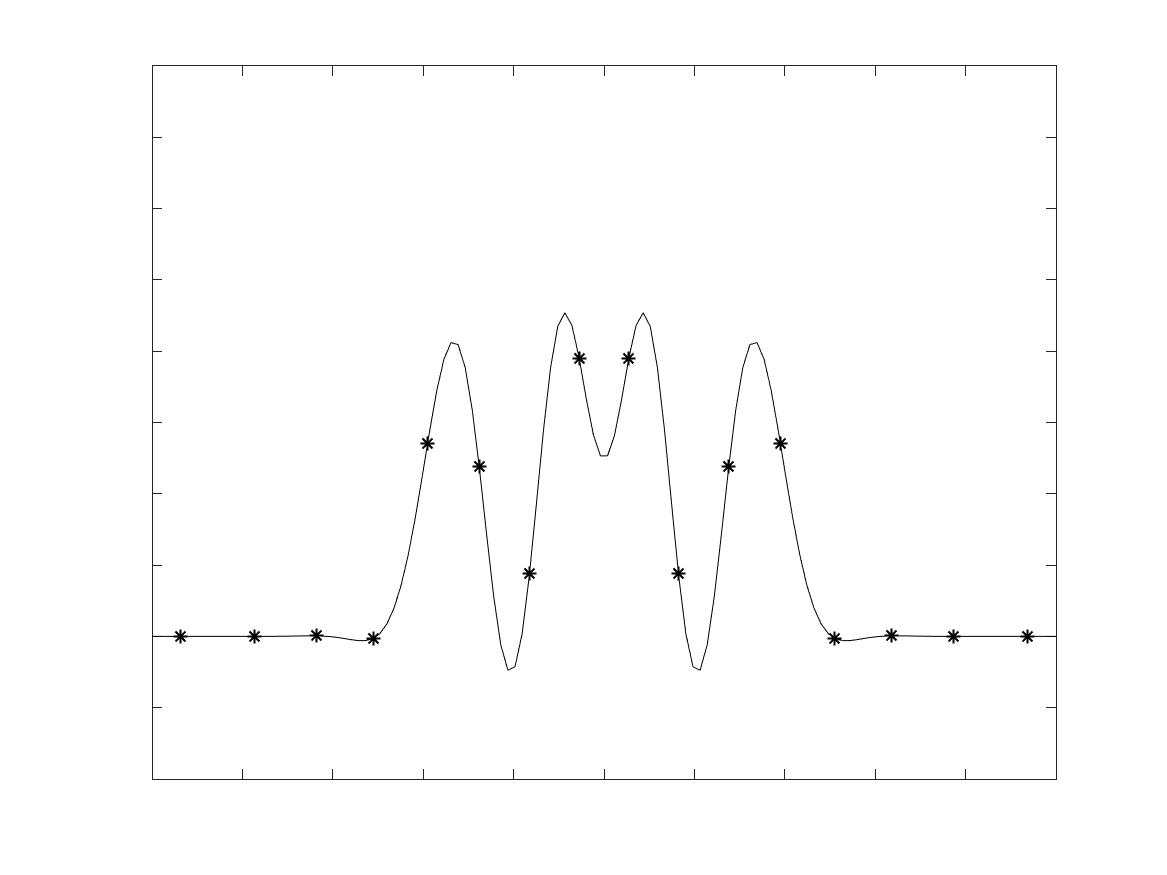}%{./figurePaper/proiezioneT30alpha1.jpg} 
  }
  \centerline{
    \includegraphics[height=5.5cm]{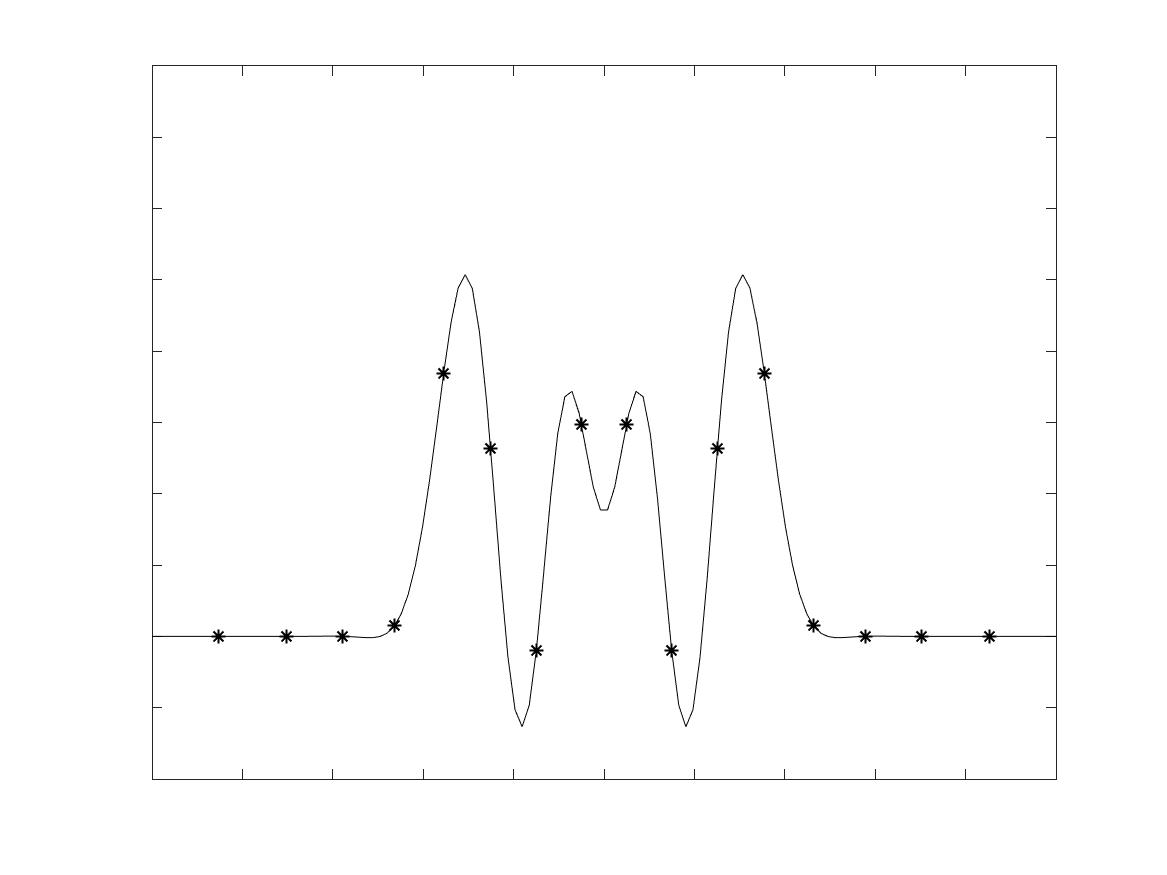}%{./figurePaper/proiezioneT30alpha11.jpg}  
    \hskip-0.6truecm    
    \includegraphics[height=5.5cm]{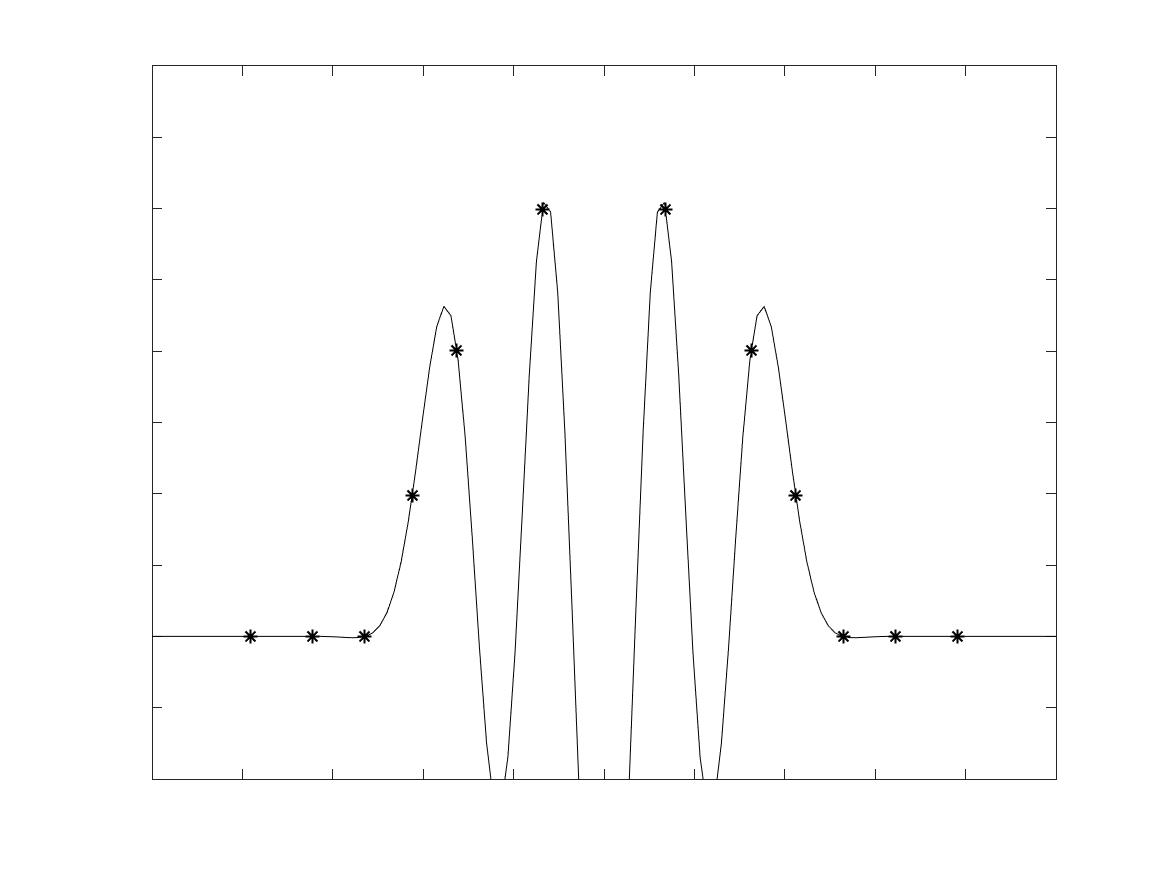}%{./figurePaper/proiezioneT30alpha12.jpg}  
  }
  \caption{ Approximated distribution function  in the variable $v$ at time $T=30$ when  $x=0$,  obtained  by using  the  Fourier-Hermite  method  with
    $N=M=2^{4}$, $\Delta t=0.01$   and    $\alpha=0.9$  (top-left), $\alpha=1$  (top-right),
    $\alpha=1.1$ (bottom-left) and $\alpha=1.2$ (bottom-right).}
  \label{fig8Sezioni}
\end{figure}

The results concerning the use of Hermite function deserve a deeper
discussion.
Figure \ref{fig7} shows the numerical distribution at time $T=25$ and
time $T=30$ obtained by using the Fourier-Hermite method with
$N=M=2^{4}$, $\Delta t=0.01$ and $\alpha=1$.
Besides, in Figure \ref{fig8} we display the approximated solution at
time $T=25$ and time $T=30$ obtained with the same parameters
(i.e. $N=M=2^{4}$, $\Delta t=0.01$) and different values of the
parameter $\alpha$.
Figure \ref{fig8Sezioni} shows the approximated distribution function
at time $T=30$ when $x=0$, obtained by using Fourier-Hermite method
with $N=M=2^{4}$, $\Delta t=0.01$ for various $\alpha$.
These figures point out that the approximated solution is very
sensitive to the choice of $\alpha$.
It has to be noted, however, that the parameter $\alpha$ that well
performs for the initial datum \eqref{eq:1D1V:twostreamf} does not
necessarily behave optimally as time increases.
For example, the value $\alpha=1.8$ looks the best choice for the
approximation of the initial distribution (see Figure \ref{fig2approx}
on bottom), nevertheless, keeping $\alpha$ constantly equal to such a
value brings to instability before arriving at time $T=30$.
This behavior clearly suggests that a dynamical way to vary $\alpha$
during the computations should be recommended, though there are at the
moment no reasonable ideas on how to realize it in practice.
%%
% These results suggests that, when scaled Hermite nodes are used, the
% best way to achieve good numerical approximations is to choose the
% value of the parameter $\alpha$ in a ``dynamic way''. Some ad hoc
% effective strategies for the parameter $\alpha$ are proposed and
% discussed in a forthcoming paper.

%% {\textcolor{red}{Aggiungere ulteriori osservazioni}}

\smallskip

\begin{figure}
   \centerline{
    \includegraphics[height=5.5cm]{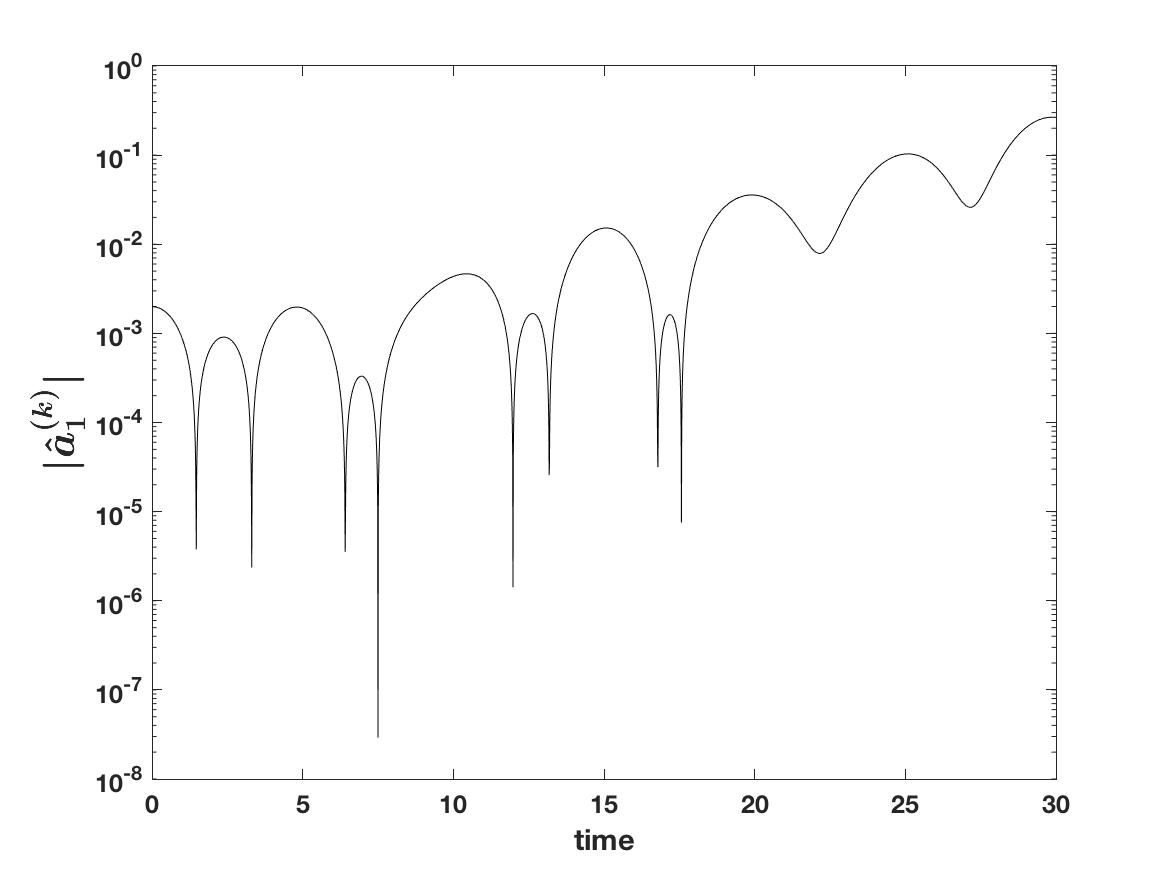}%{./figurePaper/cFFN16BDF2.jpg} 
     \hskip-0.6truecm 
    \includegraphics[height=5.5cm]{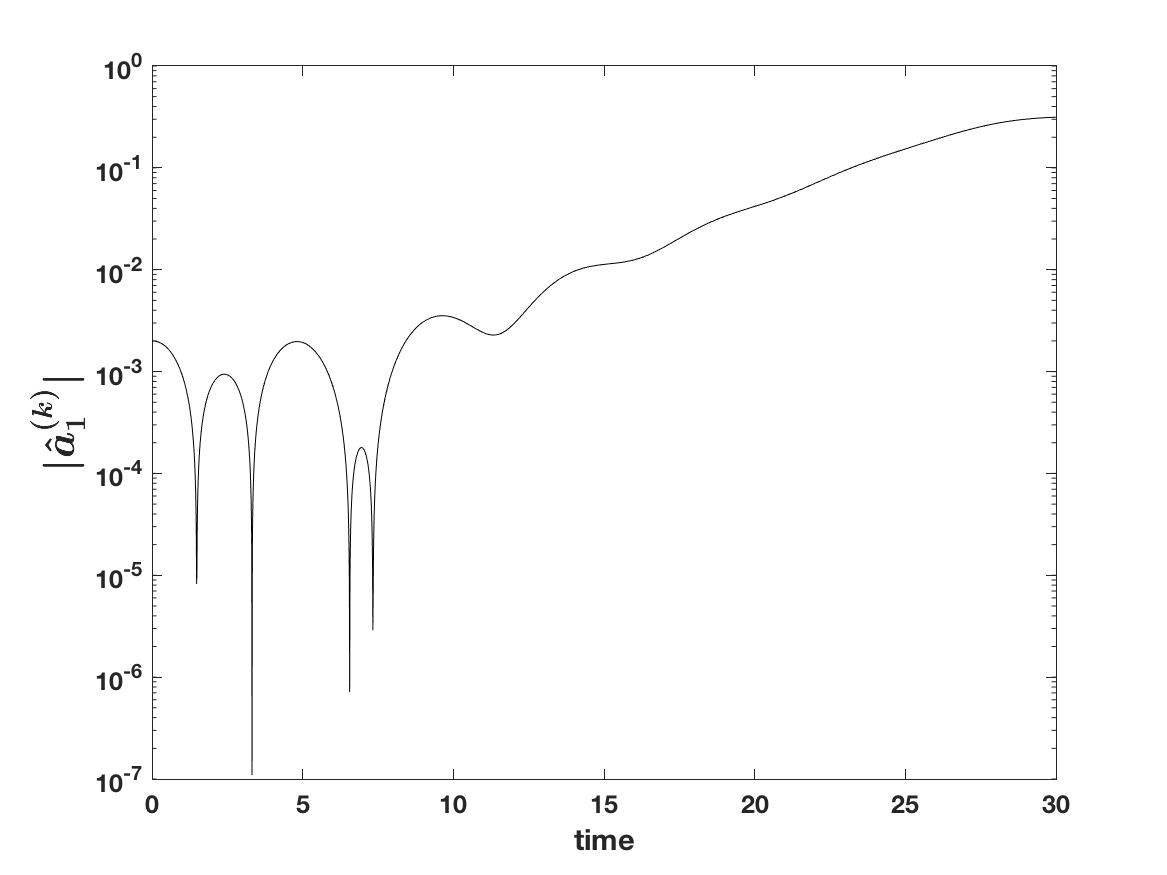}%{./figurePaper/cFHN16alpha1BDF2.jpg} 
  }
  \caption{Time evolution of the first
Fourier mode of the electric field $\Es^{(k)}_N$, i.e.
$|\hat{\as}_{1}^{(k)}|$ in \eqref{eq:1D1V:FourierSeriesE}, when using the  Fourier-Fourier method (left)  and  the  Fourier-Hermite method
(right) with    $N=M=2^{4}$, $\Delta t=0.01$ and  $\alpha=1$.}
  \label{fig9}
\end{figure}

In Figure \ref{fig9} we plot the time evolution of the (log of the)
first Fourier mode of the electric field $\Es^{(k)}_N$, i.e.
$|\hat{\as}_{1}^{(k)}|$ in ~\eqref{eq:1D1V:FourierSeriesE}, when the
Fourier-Fourier method (top) and the Fourier-Hermite method (bottom)
with $N=M=2^{4}$, $\Delta t=0.01$ and $\alpha=1$ are used.
If we compare with Figure \ref{fig4}, the two plots behave rather
differently, being the second much better that the first one.

\smallskip

\begin{figure}
  \centerline{
    \includegraphics[height=5.5cm]{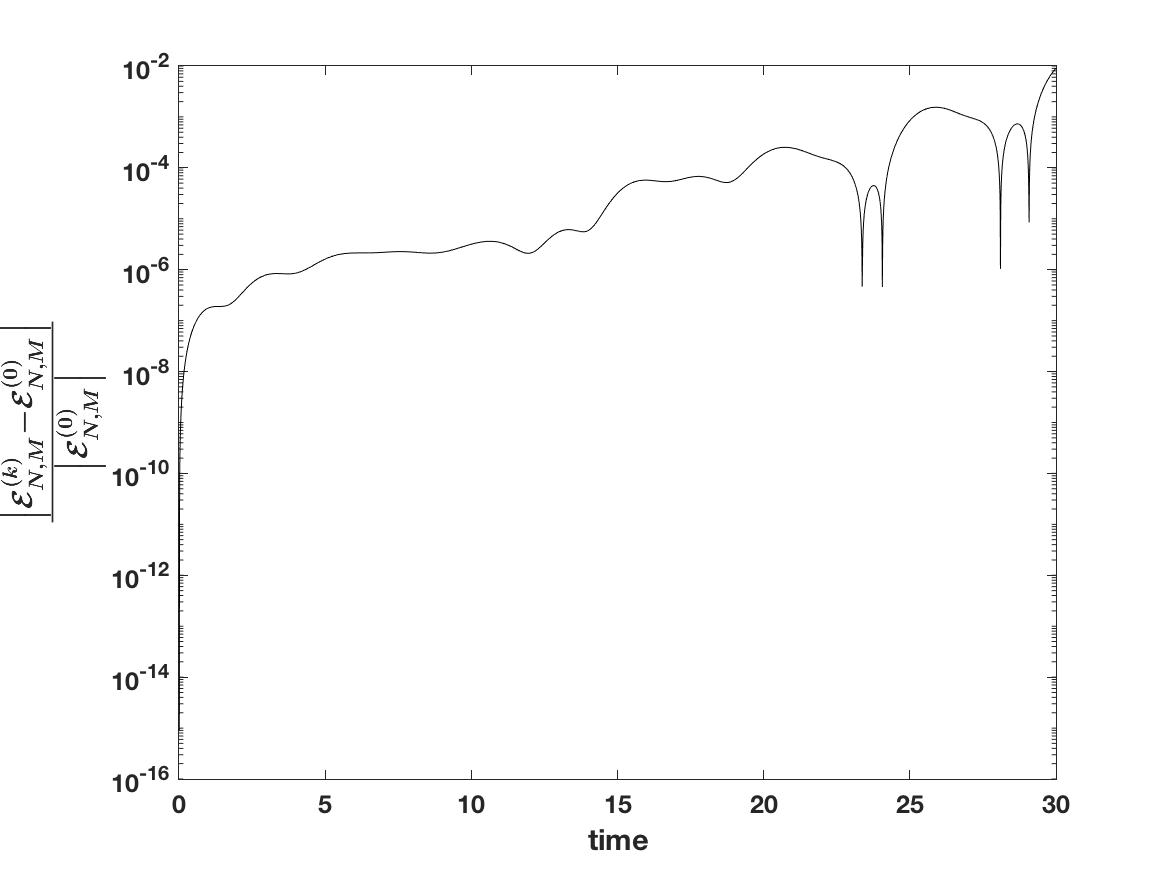}%{./figurePaper/dFFN16.jpg} 
    \hskip-0.6truecm
    \includegraphics[height=5.5cm]{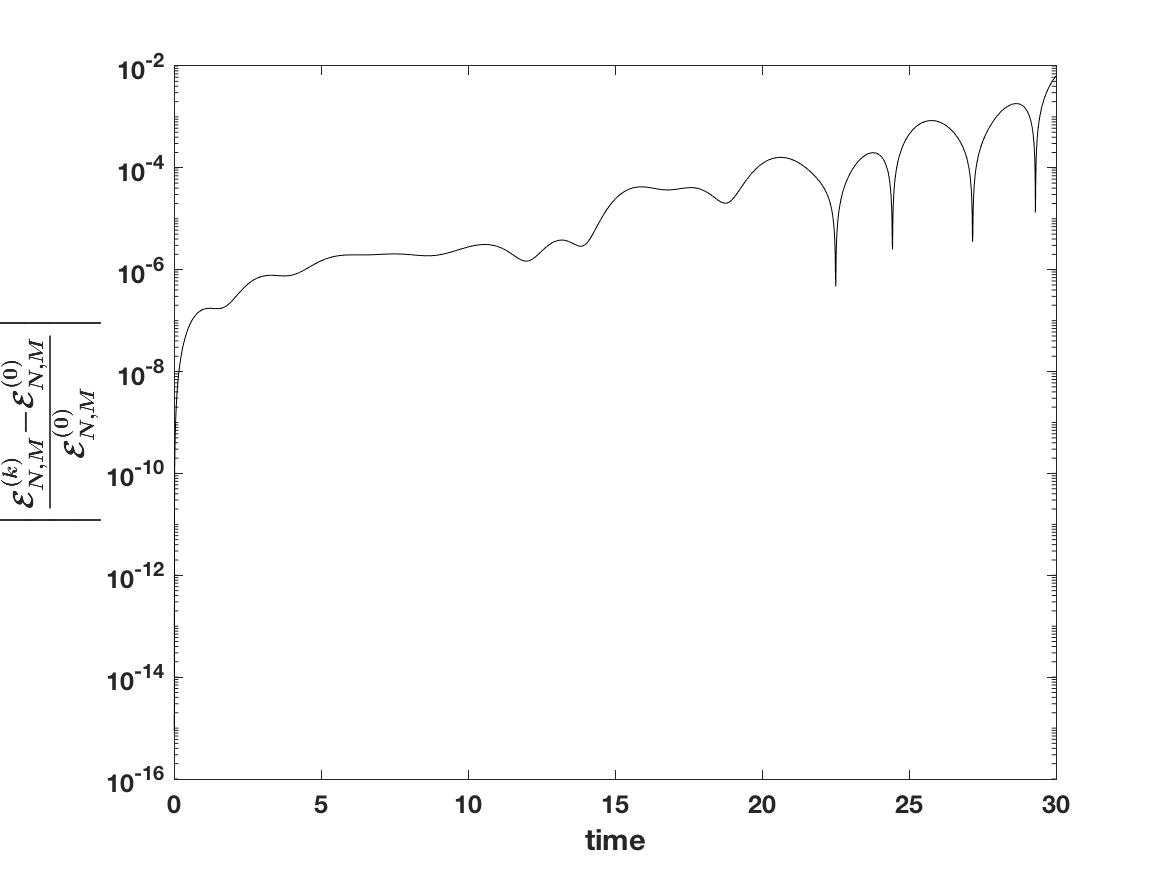}%{./figurePaper/dFFN16BDF2.jpg} 
  }
  \caption{Growth of the total energy using the first-order one-step Fourier-Fourier scheme (left)  and the second-order BDF Fourier-Fourier scheme
(right) with $N=M=2^{4}$,
and    $\Delta t=0.01$.}
  \label{fig10}
%\end{figure}
%
%\begin{figure}
  \centerline{
    \includegraphics[height=5.5cm]{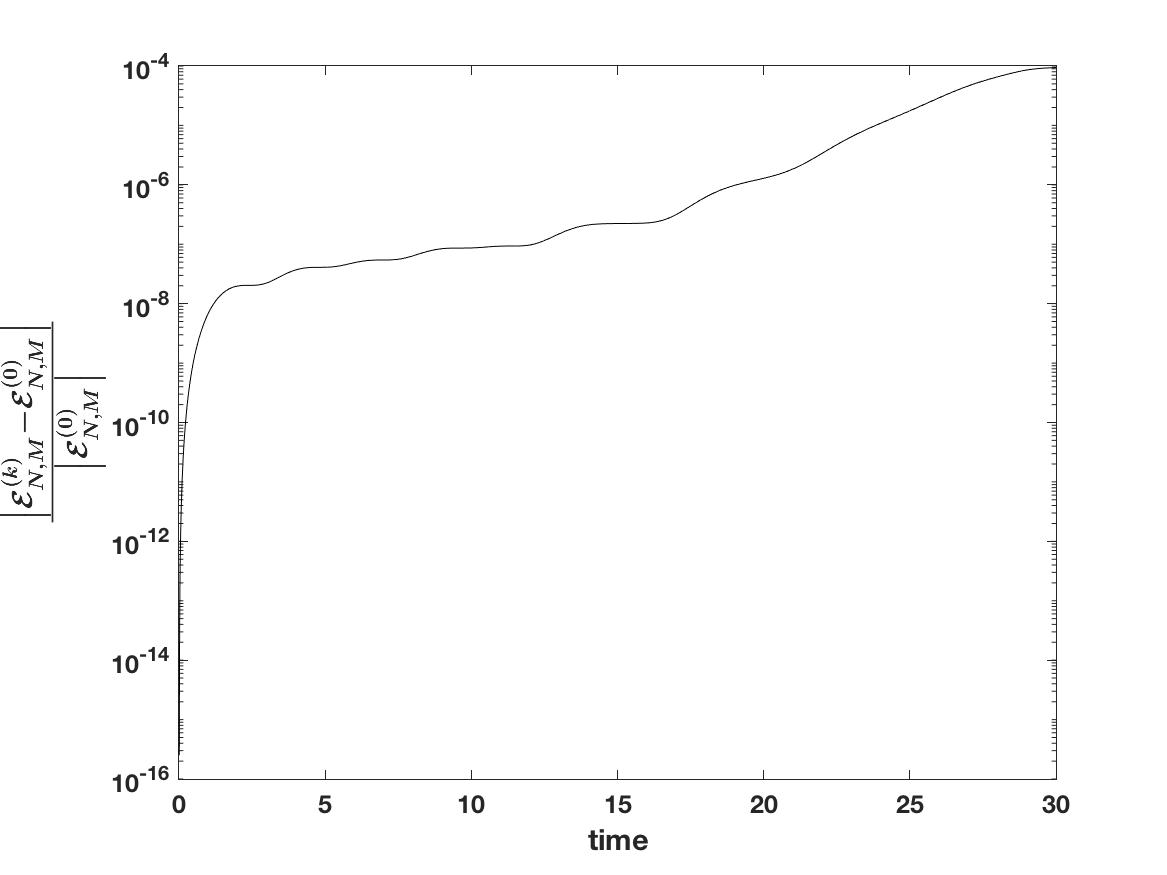}%{./figurePaper/dFHN16alpha1.jpg}    
     \hskip-0.6truecm 
    \includegraphics[height=5.5cm]{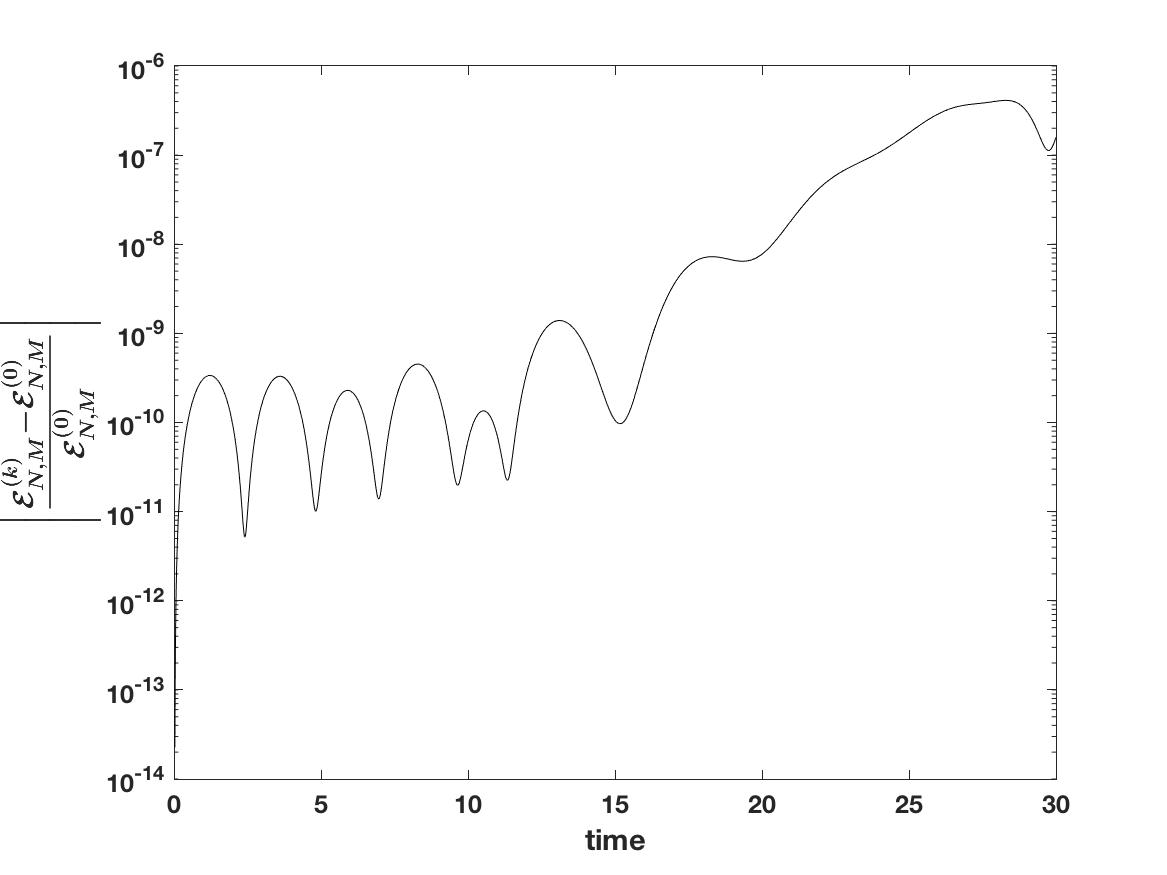}%{./figurePaper/dFHN16alpha1BDF2.jpg}       
  }
  \caption{Growth of the total energy using the first-order one-step Fourier-Hermite  scheme (left)  and the second-order BDF Fourier-Hermite scheme 
(right) with $N=M=2^{4}$,
     $\Delta t=0.01$ and  $\alpha=1$.}
  \label{fig11}
\end{figure}

In all the numerical experiments the number of particles
$\Qs^{(k)}_{N,M}$ and the momentum $\Ps^{(k)}_{N,M}$ defined in
~\eqref{eq:1D1V:massadiscreta} and ~\eqref{eq:1D1V:momentodiscreto}
are preserved.
The situation is slightly different concerning the energy, where a slight growth is observed,
with a rate depending on the size of the time step. 
Therefore, we conclude our series of tests by showing in Figures
\ref{fig10} and \ref{fig11} the variation versus time (in a
semi-$\log$ diagram) of the following quantity:
\begin{equation} \label{eq:1D1V:varEnergy}
  \frac{\left| \mathcal{\Es}_{N,M}^{(k)} - \mathcal{\Es}_{N,M}^{(0)}  \right|}{\left|\mathcal{\Es}_{N,M}^{(0)}\right|},  
\end{equation}
where the discrete energy $\mathcal{\Es}_{N,M}^{(k)}$, is defined
in~\eqref{eq:1D1V:energiadiscreta}.
Regarding energy conservation, 
for the same resolution in the phase space ($N=M=2^4$) and time
($\Delta t=0.001$), the performance of the Fourier-Hermite method
examined in this work  is superior to
that of the Fourier-Fourier method
studied in~\cite{Fatone-Funaro-Manzini:2018}.

% It is worthwhile to note that the resolution of all the numerical
% experiments shown in Figures \ref{fig5}-\ref{fig11} is $N=M=2^{4}$.
% The global number of degrees of freedom $2^4\times 2^4=256$ is
% extremely low but allows us to highlight the pros and cons of the
% various approaches.
%% 
% Noteworthy, the choice $N=M=2^{4}$ already gives reliable
% approximation results but to have sufficiently accurate results it
% is recommendable to increase $N$ up to $2^{5}$ and $M$ up to $2^{7}$
% (see~\cite{Fatone-Funaro-Manzini:2018}).
%%
%Also in this case the global number of degrees of freedom $2^5\times
%2^7=32\times128$ remains rather low.

%%%%%%%%%%%%%%%%%%%%%%%%%%%%%%%%%%%%%%%%%%%%%%%%%%%%%% 
% Soluzione di riferimento
%%%%%%%%%%%%%%%%%%%%%%%%%%%%%%%%%%%%%%%%%%%%%%%%%%%%%%

%%%%%%%%%%%%%%%%%%%%%%%%%%%%%%%%%%%%%%%%%%%%%%%%%%%%%% 
% Fourier - Fourier BDF2
%%%%%%%%%%%%%%%%%%%%%%%%%%%%%%%%%%%%%%%%%%%%%%%%%%%%%%

%%%%%%%%%%%%%%%%%%%%%%%%%%%%%%%%%%%%%%%%%%%%%%%%%%%%%% 
% Fourier -Legendre BDF2
%%%%%%%%%%%%%%%%%%%%%%%%%%%%%%%%%%%%%%%%%%%%%%%%%%%%%%

%%%%%%%%%%%%%%%%%%%%%%%%%%%%%%%%%%%%%%%%%%%%%%%%%%%%%% 
% Fourier - Hermite BDF2
%%%%%%%%%%%%%%%%%%%%%%%%%%%%%%%%%%%%%%%%%%%%%%%%%%%%%%

%%%%%%%%%%%%%%%%%%%%%%%%%%%%%%%%%%%%%%%%%%%%%%%%%%%%%% 
% Fourier - Hermite: BDF2 versus alpha
%%%%%%%%%%%%%%%%%%%%%%%%%%%%%%%%%%%%%%%%%%%%%%%%%%%%%% 

%%%%%%%%%%%%%%%%%%%%%%%%%%%%%%%%%%%%%%%%%%%%%%%%%%%%%% 
%Sezioni
%%%%%%%%%%%%%%%%%%%%%%%%%%%%%%%%%%%%%%%%%%%%%%%%%%%%%% 

%%%%%%%%%%%%%%%%%%%%%%%%%%%%%%%%%%%%%%%%%%%%%%%%%%%%%% 
% Energia
%%%%%%%%%%%%%%%%%%%%%%%%%%%%%%%%%%%%%%%%%%%%%%%%%%%%%%
\

%%%%%%%%%%%%%%%%%%%%%%%%%%%%%%%%%%%%%%%%%%%%%%%%%%%%%% 
% Energia
%%%%%%%%%%%%%%%%%%%%%%%%%%%%%%%%%%%%%%%%%%%%%%%%%%%%%% 

\section{Conclusions}
\label{sec:conclusions}

In this work, we extended the semi-Lagrangian spectral method
developed in~\cite{Fatone-Funaro-Manzini:2018},  by implementing 
Legendre polynomials and Hermite functions to approximate the distribution function with respect to the
velocity variable.
In particular, we discussed second-order accurate-in-time semi-Lagrangian methods,  obtained by
coupling spectral techniques in the space-velocity domain with a BDF time-stepping
scheme.
The resulting numerical code possesses good conservation properties,
which have been assessed by a series of numerical tests conducted
on the standard two-stream instability benchmark problem.
We also investigated numerically the dependence of the Fourier-Hermite
approximations on the scaling parameter in the Gaussian weight.
Our experiments for different representative values of this
parameter confirm that a proper choice may significantly impact on accuracy,
thus confirming previous results from the literature. The next step would be the
development of a strategy for the automatic set up of the parameter
during the time-advancing procedure.

\section*{Acknowledgements}
The second author was partially supported by the \emph{Short Term
  Mobility Program} of the Consiglio Nazionale delle Ricerche
(CNR-Italy).
The third author was supported by the Laboratory Directed Research and
Development Program (LDRD), U.S. Department of Energy Office of
Science, Office of Fusion Energy Sciences, and the DOE Office of
Science Advanced Scientific Computing Research (ASCR) Program in
Applied Mathematics Research, under the auspices of the National
Nuclear Security Administration of the U.S. Department of Energy by
Los Alamos National Laboratory, operated by Los Alamos National
Security LLC under contract DE-AC52-06NA25396.

\end{document}